\documentclass[11pt]{amsart}
\usepackage{mathtools, accents}

\usepackage{pdfpages, multirow}
\usepackage{subcaption, graphicx}
\usepackage[left=20mm,top=0.6in,bottom=12mm]{geometry}
\usepackage{enumitem}

\usepackage{xcolor}
\usepackage{hyperref}
\usepackage[hyperpageref]{backref}
\hypersetup{
    colorlinks=true,
    linkcolor=myblue,
    filecolor=magenta,      
    urlcolor=mygreen,
    citecolor=mygreen,
}
\definecolor{mygreen}{rgb}{0.01,0.5,0.2}
\definecolor{myblue}{rgb}{0.01, 0.5, 1.0}

\allowdisplaybreaks

\numberwithin{equation}{section}

\newtheorem{Lemma}{LEMMA}[section]
\newtheorem{Theorem}[Lemma]{Theorem}
\newtheorem{Proposition}[Lemma]{Proposition}
\newtheorem{Corollary}[Lemma]{Corollary}

\newtheorem{remark}[Lemma]{Remark}
\newtheorem{definition}[Lemma]{Definition}
\newtheorem{example}[Lemma]{Example}

\newtheorem{Fact}[Lemma]{Fact}
\newtheorem{assumption}[Lemma]{Assumption}

\def\bt{\begin{Theorem}}
\def\et{\end{Theorem}}
\def\bl{\begin{Lemma}}
\def\el{\end{Lemma}}
\def\bp{\begin{Proposition}}
\def\ep{\end{Proposition}}
\def\bcor{\begin{Corollary}}
\def\ecor{\end{Corollary}}
\def\bpf{\begin{proof}}
\def\epf{\end{proof}}

\def\brem{\begin{remark}\rm }
\def\erem{\hfill $\lozenge$ \end{remark}}

\def\bedef{\begin{definition}\rm }
\def\endef{\hfill$\lozenge$\end{definition}}

\def\beg{\begin{example}\rm }
\def\eeg{\hfill $\lozenge$\end{example}}

\def\bef{\begin{Fact}}
\def\eef{\end{Fact}}

\def\bea{\begin{assumption}}
\def\ena{\end{assumption}}
\def\bc{\begin{center}}
\def\ec{\end{center}}
\def\noi{\noindent}

\def\beq{\begin{equation}}
\def\eeq{\end{equation}}
\def\beqarray{\begin{eqnarray*}}
\def\eeqarray{\end{eqnarray*}}
\def\<{\leftangle}
\def\>{\rightangle}
\def\({\left(}
\def\){\right)}
\def\f{\varphi}

\def\<{\langle}
\def\>{\rangle}
\def\q{\quad}

\def\a{\alpha}
\def\b{\beta}
\def\g{\gamma}

\def\d{\delta}

\def\k{\kappa}

\def\t{\tau}

\def\l{\lambda}

\def\O{\Omega}
\def\o{\omega}

\def\w.r.t.{with respect to}
\def\R{{\mathbb{R}}}
\def\N{{\mathbb{N}}}

\def\C{{\mathbb{C}}}

\def\bq{\begin{quote}}
\def\eq{\end{quote}}

\def\bit{\begin{itemize}}
\def\eit{\end{itemize}}
\def\i{\item}
\def\ben{\begin{enumerate}}
\def\een{\end{enumerate}}

\def\Ea{{E_{\a,1}}}
\def\Eaa{{E_{\a,\a}}}
\def\Egb{{ E_{\g,\b}}}

\def\dat{{\partial^\a_t}}

\begin{document}
\title[Backward problem for time-fractional fourth order equation]{On backward problem for a time-fractional fourth order parabolic equation}
\author{Subhankar Mondal}
%
\address{(SM): TIFR Centre for Applicable Mathematics, Bangalore-560065, India;}
\email{subhankar22@tifrbng.res.in; }

\begin{abstract} 
This article is concerned with the inverse problem of retrieving the initial value of a time-fractional fourth order parabolic equation from source and final time observation. The considered problem is an {\it ill-posed problem.} We obtain regularized approximations for the sought initial value by employing the quasi-boundary value method, its modified version and by Fourier truncation method(FTM). We provide both the apriori and aposteriori parameter choice strategies and derive the error estimates for all these methods under some {\it source conditions} involving some Sobolev smoothness. As an important implication of the obtained rates, we observe that for both the apriori and aposteriori cases, the rates obtained by all these three methods are same for some source sets. Moreover, we observe that in both the apriori and aposteriori cases, the FTM is free from the so-called {\it saturation effect}, whereas the quasi-boundary value method and its modifications possesses the saturation effect for both the cases. Further, we observe that the rates obtained by the FTM is always order optimal for the considered source sets.  
\end{abstract}

\maketitle

\noi\textbf{Keywords:} backward problem, time-fractional parabolic equation, regularization, quasi boundary value method, truncation method, convergence rate

\noi\textbf{MSC 2020:} 47A52, 35R25, 35R11, 65J20
\section{Introduction} 
Fourth order parabolic PDE's have remained an interesting topic of research for many decades. It is known that fourth order parabolic PDE's are helpful in modelling various natural phenomenon and has found its application in many engineering problems such as, the free vibrations in beams or shafts \cite{gomain_1975}, the growth of epitaxial thin film \cite{king_stein_winkler_2010}, thermal grooving by surface mechanism \cite{mullins_1957}, and is also used for noise removal and edge preservation in image processing \cite{you_kaveh_2000}, to name a few. Recently research activities on inverse problems related to fourth order parabolic PDE has grown up substantially. Let us mention briefly a few recent works related to some inverse problems associated with fourth order parabolic equations. 

For $\t>0,$ and some appropriate functions $h,v_0,$ in \cite{cao_lesnic_ismailov_2021} the authors considered the system of equations
\begin{equation*}
\begin{cases}
v_t(x,t)+a(t) v_{xxxx}(x,t) = h(x,t),&\q \text{in}\q (0,1)\times(0,\t),\\
v_x(0,t)=v_x(1,t)=v_{xxx}(0,t)=v_{xxx}(1,t)=0,&\q \text{on}\q (0,\t),\\
v(x,0)=v_0(x),&\q \text{in}\q (0,1)
\end{cases}
\end{equation*}
and analysed the inverse problem of identifying the time-dependent coefficient $a(t),$ called the {\it Mullins' coefficient,} from the knowledge of the data $v(x_0,t),$ where $x_0\in [0,1]$ is any fixed given point and $v$ is the solution of the considered PDE. Subsequently, in \cite{cao_AMO_2022}, the author considered a more general fourth order parabolic PDE and investigated the inverse source identification problem associated with it. More precisely, the author considered the inverse problem of identification of the source function $f$ from the knowledge of final time observation $v(x,\t)$ or the time-integral observation $\int_0^\t \o(t)v(x,t)\,dt,$ for $0<x<1,$ where $v$ is the solution of 
\begin{equation*}
\begin{cases}
v_t(x,t)+v_{xxxx}(x,t) = f(x) g(x,t)+ h(x,t),&\q \text{in}\q (0,1)\times(0,\t),\\
v(0,t)=\mu_1(t),\,\,v(1,t)=\mu_2(t),&\q \text{on}\q (0,\t),\\
v_{xx}(0,t)=\mu_3(t),\,\,v_{xx}(1,t)=\mu_4(t),&\q \text{on}\q (0,\t),\\
v(x,0)=v_0(x),&\q \text{in}\q (0,1),
\end{cases}
\end{equation*}
for some given functions $\o,g,h,v_0,\mu_i,i=1,2,3,4.$

As a follow up to the above developments, in a recent work \cite{huntul_2024}, the authors have extended the investigation from integer-order derivative to fractional order setting. They considered an inverse problem of identifying time-dependent source term $f$ from the final value $v(\cdot,\t),$ associated with a time-fractional fourth order parabolic PDE as follows:
\begin{equation*}
\begin{cases}
\dat v(x,t)+B v_{xxxx}(x,t) = f(t) g(x,t),&\q \text{in}\q (0,1)\times(0,\t),\\
v(1,t)=v_x(0,t)=v_{xx}(1,t)=v_{xxx}(0,t)=0,&\q \text{on}\q (0,\t),\\
v(x,0)=v_0(x),&\q \text{in}\q (0,1),
\end{cases}
\end{equation*}
where $B>0$ is a surface diffusion coefficient, $v_0$ is some appropriate function. Here (and throughout the paper) $\dat$ is the Caputo fractional derivative of order $\a$ (cf. \cite{podlubny_1999_book}) given by
\beq\label{caputo_def}
\dat v(x,t):=\begin{cases}
\frac{1}{\Gamma(1-\a)}\int_0^t \frac{v_s(x,s)}{(t-s)^\a}\,ds,\q&\text{if}\q 0<\a<1,\\
v_t(x,t),\q&\text{if}\q \a=1,
\end{cases}
\eeq
where $\Gamma(\cdot)$ is the standard Gamma function, that is, $\Gamma(t)=\int_0^\infty e^{-s}s^{t-1}\,ds,\,\,t>0.$

Note that forward problem associated with such time-fractional fourth order parabolic equations have been studied by many authors in various context, see for instance \cite{hamed_nepo_2015, hristov_2018}. 

It is important to recall that over the last few decades the research on time-fractional backward heat conduction problems (TFBHCP) have gained substantial attention due to its many practical applications. Since the work in \cite{liu_yamamoto_2010}, the literature on TFBHCP and its various variants have grown remarkably (see e.g., \cite{hao_liu_duc_thang_2019,kokila_nair_2020, tuan_long_nguyen_tran_2017}). Considering the importance of fractional-order derivative, it is desirable to study various inverse problem associated with time-fractional fourth order parabolic equations. Recently some works on source identification problem associated with time-fractional fourth order parabolic equations have been reported (see e.g. \cite{huntul_2024}), as mentioned earlier. However, to the best of our knowledge, the backward problem for time-fractional fourth order parabolic equations has not been explored yet and thus, this paper intends to contribute in that direction.

In this paper we are interested in the inverse problem of retrieving the initial value from the knowledge of final value and the source function, associated with a time-fractional fourth order parabolic PDE. More precisely, we consider the following: Let $\O\subset \R^d,d\in \{1,2,\ldots,7\},$ (see Remark \ref{dimension-restriction} for this restriction on $d$) be a bounded domain with sufficiently smooth boundary $\partial\O$ and $\t>0$ be fixed. Let $\a\in (0,1).$ We consider the system
\beq\label{gov-pde}
\begin{cases}
\partial^\a_t u+\Delta^2 u = f,&\q\text{in}\q\O\times (0,\t),\\
u=0=\Delta u,&\q \text{on}\q\partial\O\times (0,\t),\\
u(\cdot,\t)=h,&\q\text{in}\q\O,
\end{cases}
\eeq
where $f\in L^\infty([0,\t];L^2(\O)),$ $h\in L^2(\O)$ and $\dat$ is the Caputo fractional derivative as defined in \eqref{caputo_def}. The inverse problem that we consider is:
\bq
{\bf (IP)}: retrieve the initial value $u(\cdot,0)$ from the knowledge of $h$ and $f,$ where $u$ is a solution of \eqref{gov-pde}.
\eq  

Similar to TFBHCP, the inverse problem to recover $u(\cdot,0)$ from the knowledge of $h,f,$ associated with \eqref{gov-pde} is an ill-posed problem in the sense that a small perturbations in the data may lead to a large deviation in the sought initial value (see Example \ref{eg-illposed}). Further, the problem of recovering $u(\cdot,t)$ for $0<t<\t$ is stable \w.r.t. the perturbations in the data $h,f,$ (see Remark \ref{rem-stability-for-t-non-zero}), possessing the properties similar to that of TFBHCP.

In practice, we shall have access to some measured data, and hence noise in the data is inevitable. We assume that for $\d>0,$ $h^\d\in L^2(\O)$ and $f^\d\in L^\infty(0,\t;L^2(\O))$ are the noisy data satisfying 
\beq\label{noise-model}
\|h-h^\d\|_{L^2(\O)}^2+\theta\|f-f^\d\|^2_{L^\infty(0,\t;L^2(\O))}\leq \d^2,
\eeq
 where $\theta>0$ is a fixed known constant that will be specified later (see Proposition \ref{combined-exact-data}). 

Thus, the inverse problem at hand is to recover $u(\cdot,0)$ from the knowledge of $h^\d$ and $f^\d$ satisfying \eqref{noise-model}. Since the problem is ill-posed, in order to obtain stable approximations, we have to employ some regularization methods (see e.g., \cite{engl_hanke_neubauer,nair_opeq}). We obtain stable approximations by the quasi-boundary value method (QBVM) (see e.g., \cite{hao_duc_lesnic_2010, tuan_long_nguyen_tran_2017}), it's modifications (MQBVM) (see e.g., \cite{hao_liu_duc_thang_2019, mondal_2024, wei_wang_2014}) and the Fourier truncation method (FTM) (see e.g., \cite{kokila_nair_2020, mondal_JOTA_2024}). Just to give some insight, in QBVM we perturb the final value in \eqref{gov-pde} by adding the initial value multiplied by the regularization parameter, that is, we consider the resulting perturbed term as $u(\cdot,\t)+\b u(\cdot,0),$ where $\b>0$ is the regularization parameter. For the MQBVM  we consider the resulting perturbed term as $u(\cdot,\t)+\b (-\Delta)^q u(\cdot,0),$ where $\b>0$ is the regularization parameter and $q\in\N$. Note that in literature, the case $q=1,$ has been studied extensively for various problems (see e.g., \cite{wei_wang_2014}), and is known by the name  modified quasi-boundary value method. In FTM, we consider the truncated version of the Fourier expansion of the concerned function and the level at which we truncate plays the role of regularization parameter.

It is known in regularization theory that unless we assume some {\it source condition} (see e.g., \cite{engl_hanke_neubauer, nair_opeq}) on the unknown, the convergence can be arbitrarily slow \cite{schock_1984}, that is, we can not obtain any convergence rate. In view of  this, we assume that the unknown $g=u(\cdot,0)$ belongs to some {\it source sets} which are subsets of some Sobolev spaces $H_p,$ where $p>0$ indicates the smoothness property, and then obtain error estimates and hence deduce the convergence rates of the approximations by the considered regularization methods. We shall provide both the apriori and the aposteriori parameter choice strategies and obtain the corresponding error estimates for all of them. As an important implication, we observe that for the FTM, the convergence rates in both the apriori and aposteriori cases is of the form $O(\d^{p/(p+2)})$ for all $p>0,$ where $\d$ is the noise as considered in \eqref{noise-model}. However, for the QBVM and MQBVM, for both the apriori and aposteriori cases, the convergence rates possesses saturation effect. For the QBVM the best rates possible for the apriori and aposteriori cases is $O(\d^{1/2})$ and for MQBVM the best rates possible for the apriori and aposteriori cases are $O(\d^{(q+2)/(q+4)})$ and $O(\d^{q/(q+2)}),$ respectively, for $q\in\N.$ Further, we observe that, if $p$ belongs to certain interval then the rates of FTM, QBVM, MQBVM are all the same for apriori case, and for certain values of $p$ the aposteriori rates for FTM and MQBVM are also the same. Moreover, we shall see that for the considered source sets the rates obtained by the FTM is always order optimal for all $p>0,$ but the same is not the case for QBVM and MQBVM (see Remark \ref{optimality}).

The remainder of the paper is organised as follows: In Section \ref{sec-prelim} we list some notations, important facts, results about the {\it Mittag-Leffler function} and some basic result on calculus, that will be used throughout. In Section \ref{sec-prep-results} we mainly discuss the relation between the sought initial value and the data, that is, the final value and the source function. We show by an example the ill-posedness of the considered inverse problem of retrieving the initial value. We also establish a conditional stability result. Section \ref{sec-quasi-bound} is all about the QBVM and MQBVM. We provide both apriori and aposteriori parameter choice strategies and obtain the corresponding error estimates. Section \ref{sec-FTM} is about the FTM. We systematically define the regularized solutions and provide both apriori and aposteriori parameter choice strategies and obtain the corresponding convergence rates. Moreover, we discuss about the optimality of the obtained rate.

\section{Preliminaries}\label{sec-prelim}
In this section we recall and record some of the results, facts, notations that will be used throughout the paper. To begin with, we first recall a standard result from spectral theory of elliptic equations. Let $\O$ be as earlier. Then  the following eigenvalue problem
\beq\label{eig-prob}
\begin{cases}{}
-\Delta \f = \l\f &\q\text{in}\q\O,\\
\q\q\f =  0 &\q\text{on}\q\partial\O.
\end{cases}
\eeq
admits a sequence $(\lambda_n)$ of eigenvalues such that  $0<\l_1< \l_2\leq\ldots\to\infty$ and the corresponding sequence $(\f_n)$ of eigenfunctions  are such that $\f_n\in H^1_0(\O)\cap H^2(\O)$ which form an orthonormal basis of $L^2(\O)$ (cf. \cite{evans}). Moreover, it is known that (see e.g., \cite{courant_hilbert_book}) there exist constants $e_1,e_2>0$ such that
\beq\label{eig-value-bound}
e_1n^{2/d}\leq \l_n\leq e_2 n^{2/d},\q\forall\,n\in\N.
\eeq

Before proceeding further, let us list the following notations that we will be followed throughout.
\bit 
\i For $z\in L^2(0,\t;L^2(\O)),$ we may use the notation $z(t)$ to denote $ z(\cdot,t)$ for almost all $t\in [0,\t].$
\i $C([0,\t];L^2(\O))$ denotes the standard space of all $L^2(\O)$-valued continuous function in $[0,\t].$
\i For $g_1,g_2\in L^2(\O)$, $\<g_1,g_2\>$ denotes the standard inner-product in $L^2(\O)$ and we shall denote the norm $\|\cdot\|_{L^2(\O)}$ by $\|\cdot\|.$ 
\i For $\phi\in L^2(\O),$ we shall always denote $\<\phi,\f_n\>$ by $\phi_n$ for all $n\in\N.$
\i For $z\in L^\infty(0,\t;L^2(\O))$ we shall denote $\<z(s),\f_n\>$ by $z_n(s)$ for $0\leq s\leq \t$ and for all $n\in\N.$
\eit

We now recall the definition of Mittag-Leffler function (cf. \cite{podlubny_1999_book}). For $\g>0$ and $\b\in \R,$ the Mittag-Leffler function denoted by $\Egb,$ is defined by
$$\Egb(z):=\sum_{k=0}^\infty \frac{z^k}{\Gamma(k\g+\b)},\q z\in\C,\,\g>0,\,\b\in\R.$$
\bt\label{prop_Ea_upper_bound_lower_bound}\cite{liu_yamamoto_2010}
Let $\a\in (0,1).$ There exist constants $C_1,C_2>0$ such that
$$\frac{C_1}{\Gamma(1-\a)}\frac{1}{(1-x)}\leq \Ea(x)\leq \frac{\bar{C}_1}{\Gamma(1-\a)}\frac{1}{(1-x)}\q\forall x\leq 0.$$
\et
\bt\label{prop_Eaa_L1bound}\cite[(3.7), pg. 434]{sakamoto_yamamoto_2011}
Let $\a\in (0,1)$ and $\l>0.$ The function $t\mapsto t^{\a-1}\Eaa(-\l t^\a),$ belongs to $L^1(0,\t)$ and $$\int_0^\t |t^{\a-1}\Eaa(-\l t^\a)|\,dt\leq \frac{1}{\l}.$$
\et
We end this section by stating a result whose proof follows from basic calculus.
\bl\label{calculus}
Let $c>0$ be fixed. The following holds:
\ben
\i[(i)] For $\b,q>0$, the function $\psi_{q,\b}:[\l_1,\infty)\to [0,\infty)$, defined by $$\psi_{q,\b}(s)=\frac{s^2}{c+\b s^{q+2}}$$ attains its global maximum at $s_0=(2c/q)^{1/(q+2)}\frac{1}{\b^{1/(q+2)}}$ and $\psi(s_0)=(2c/q)^{2/(q+2)}\frac{q}{c(1+q)}\frac{1}{\b^{2/(q+2)}}.$
\i[(ii)] Let $q\geq 0$ and $\b>0.$ For $0<p<q+2,$ the function $\phi_{p,q,\b}:[\l_1,\infty)\to [0,\infty)$ defined by $$\phi_{p,q,\b}(s)=\frac{\b s^{q+2-p}}{c+\b s^{q+2}}$$ attains its global maximum at $s_0=(c(q+2-p)/p)^{1/(q+2)}\frac{1}{\b^{1/(q+2)}}$ and $$\phi_{p,q,\b}(s_0)=(c(q+2-p)/p)^{(q+2-p)/(q+2)}\frac{p}{c(q+2)}\b^{p/(q+2)}.$$
\een
\el
\section{Preparatory results}\label{sec-prep-results}
In order to establish a relation between the sought initial value and the data, that is, the source function and the final value, we first consider the solution representation of a forward problem. We consider the PDE
\beq\label{forward-pde}
\begin{cases}
\partial^\a_t v+\Delta^2 v = f,&\q\text{in}\q\O\times (0,\t),\\
v=0=\Delta v,&\q \text{on}\q\partial\O\times (0,\t),\\
v(\cdot,0)=g,&\q\text{in}\q\O,
\end{cases}
\eeq
where $f\in L^\infty(0,\t;L^2(\O))$ and $g\in L^2(\O).$

Let $v$ satisfies \eqref{forward-pde}. Then following \cite{huntul_2024}, we have
\beq\label{forward-sol-rep}
v(x,t)=\sum_{n=1}^\infty \big[\Ea(-\l_n^2 t^\a)g_n+\int_0^t (t-s)^{\a-1}\Eaa(-\l_n^2(t-s)^\a)f_n(s)\,ds\big]\f_n(x),
\eeq
where recall that $g_n=\<g,\f_n\>$ and $f_n(s)=\<f(\cdot,s),\f_n\>.$
\bedef
Any $v\in C([0,\t];L^2(\O))$ satisfying \eqref{forward-sol-rep} is said to be a solution of \eqref{forward-pde}.
\endef
Clearly, $v$ given by \eqref{forward-sol-rep} is indeed a solution of \eqref{forward-pde}. Now suppose that it is given $v(\cdot,\t)=h(\cdot).$ Then, we obtain
$$h(x)=v(x,\t)=\sum_{n=1}^\infty \big[\Ea(-\l_n^2 \t^\a)g_n+\int_0^\t (\t-s)^{\a-1}\Eaa(-\l_n^2(\t-s)^\a)f_n(s)\,ds\big]\f_n(x),$$ which implies 
\beq\label{g-h-f-relation-aux}
g_n=\frac{h_n-\int_0^\t (\t-s)^{\a-1}\Eaa(-\l_n^2(\t-s)^\a)f_n(s)\,ds}{\Ea(-\l_n^2\t^\a)},\,\,\forall n\in\N.
\eeq In order to reduce the notational complexity, for any $h,h^\d\in L^2(\O)$ and $f,f^\d\in L^\infty(0,\t;L^2(\O)),$ we define the following
\beq\label{psifan}
\begin{cases}
&\Psi_f^{\a,n}(t):=\int_0^t (t-s)^{\a-1}\Eaa(-\l_n^2(t-s)^\a)f_n(s)\,ds,\q 0\leq t\leq \t,\\
&\Psi_{f^\d}^{\a,n}(t) := \int_0^t (t-s)^{\a-1}\Eaa(-\l_n^2(t-s)^\a)f^\d_n(s)\,ds,\q 0\leq t\leq \t,\\
&\Upsilon_{h,f}^{\a,n}:= h_n-\Psi^{\a,n}_{f}(\t),\\
& \Upsilon_{h^\d,f^\d}^{\a,n}:=h^\d_n-\Psi^{\a,n}_{f^\d}(\t),
\end{cases}
\eeq
where we recall that $h_n=\<h,\f_n\>$ and $f_n(s)=\<f(s),\f_n\>.$

Thus, from \eqref{g-h-f-relation-aux}, we have
\beq\label{g-h-f-relation}
g_n=\frac{h_n-\Psi_f^{\a,n}(\t)}{\Ea(-\l_n^2\t^\a)}=\frac{\Upsilon_{h,f}^{\a,n}}{\Ea(-\l_n^2\t^\a)},
\eeq
and hence from \eqref{forward-sol-rep}, we have
\beq\label{sol-rep-bacward-all-time-aux}
v(x,t)=\sum_{n=1}^\infty \big[\frac{\Ea(-\l_n^2t^\a)}{\Ea(-\l_n^2\t^\a)}(h_n-\Psi^{\a,n}_f(\t))+\Psi^{\a,n}_f(t)\big]\f_n(x),\q0\leq t\leq\t.
\eeq
\bedef\label{def-sol-backward-problem}
Let $h\in L^2(\O)$ and $f\in L^\infty(0,\t;L^2(\O)).$ Any $v\in C([0,\t];L^2(\O))$ that satisfies \eqref{sol-rep-bacward-all-time-aux} is said to be a solution of the backward problem \eqref{gov-pde}.
\endef

Before proceeding further, let us put in place some technical results. 
\bp\label{psi-alphan-difference-estimate} 
Let $f,\tilde{f}\in L^\infty(0,\t;L^2(\O)).$ Let $\Psi_{f}^{\a,n}(t)$ and $\Psi_{\tilde{f}}^{\a,n}(t)$ be as defined in \eqref{psifan}. Then the following holds:
\ben
\i[(i)] For all $0\leq t\leq \t,$ we have $|\Psi_{f}^{\a,n}(t)|\leq \frac{1}{\l_n^2}\|f\|_{L^\infty(0,\t;L^2(\O))}.$
\i[(ii)] For all $0\leq t\leq \t,$ we have $|\Psi^{\a,n}_f(t)-\Psi^{\a,n}_{\tilde{f}}(t)|\leq \frac{1}{\l_n^2}\|f-\tilde{f}\|_{L^\infty(0,\t;L^2(\O))}$
\een
\ep
\bpf
{\bf (i).} We have 
\beqarray
|\Psi_{f}^{\a,n}(t)|&=&|\int_0^t (t-s)^{\a-1}\Eaa(-\l_n^2(t-s)^\a)f_n(s)\,ds|\\
&\leq &\|f\|_{L^\infty(0,\t;L^2(\O))}\int_0^\t|t^{\a-1}\Eaa(-\l_n^2 t^\a)|\,dt\\
&\underbrace{\leq}_{Theorem \ref{prop_Eaa_L1bound}} & \frac{1}{\l_n^2} \|f\|_{L^\infty(0,\t;L^2(\O))}.
\eeqarray
{\bf (ii).} The proof follows from (i).
\epf
\bp\label{combined-exact-data}
Let $\theta=\sum_{n=1}^\infty \frac{1}{\l_n^4}.$ Let $h,\tilde{h}\in L^2(\O)$ and $f,\tilde{f}\in L^\infty(0,\t;L^2(\O)).$ Let $\Upsilon_{h,f}=\sum_{n=1}^\infty \Upsilon_{h,f}^{\a,n}\,\f_n$ and $\Upsilon_{\tilde{h},\tilde{f}}=\sum_{n=1}^\infty \Upsilon_{\tilde{h},\tilde{f}}^{\a,n}\,\f_n,$ where $\Upsilon_{h,f}^{\a,n}$ is as defined in \eqref{psifan}. Then the following holds:
\ben
\i[(i)] We have $\Upsilon_{h,f}\in L^2(\O)$ and $$\|\Upsilon_{h,f}\|^2\leq 2(\|h\|^2+\theta\|f\|^2_{L^\infty(0,\t;L^2(\O))}).$$
\i[(ii)] We have $$\|\Upsilon_{\tilde{h},\tilde{f}}-\Upsilon_{h,f}\|^2\leq 2(\|h-\tilde{h}\|^2+\theta\|\tilde{f}-f\|^2_{L^\infty(0,\t;L^2(\O))}).$$
\een
\ep
\bpf
First we observe that 
\beqarray
\sum_{n=1}^\infty|\Upsilon_{h,f}^{\a,n}|^2&=&\sum_{n=1}^\infty |h_n-\Psi_{f}^{\a,n}(\t)|^2\leq 2\sum_{n=1}^\infty (|h_n|^2+|\Psi_{f}^{\a,n}(\t)|^2)\\
&\underbrace{\leq }_{Proposition \ref{psi-alphan-difference-estimate}}& 2(\|h\|^2+\theta\|f\|^2_{L^\infty(0,\t;L^2(\O))}).
\eeqarray
Now the proof follows.
\epf
\brem\label{dimension-restriction}
From \eqref{eig-value-bound}, we know that $\frac{1}{\l_n}\sim \frac{1}{n^{2/d}}.$ Thus, if $1\leq d\leq 7,$ then $\sum_{n=1}^\infty\frac{1}{\l_n^4}<\infty$ and hence it make sense to define $$\theta=\sum_{n=1}^\infty \frac{1}{\l_n^4}.$$
\erem

We are now in a position to state and prove the necessary and sufficient condition for the existence and uniqueness of the solution of backward problem \eqref{gov-pde}.
\bt\label{nec-suff-sond-backward-problem}
Let $h\in L^2(\O)$ and $f\in L^\infty(0,\t;L^2(\O)).$ The backward problem associated with \eqref{gov-pde} has a unique solution if and only if $$\sum_{n=1}^\infty \left(\frac{h_n-\Psi_{f}^{\a,n}(\t)}{\Ea(-\l_n^2\t^\a)}\right)^2<\infty$$ and the solution is given by
\beq\label{sol-rep-bacward-all-time}
u(x,t)=\sum_{n=1}^\infty \big[\frac{\Ea(-\l_n^2t^\a)}{\Ea(-\l_n^2\t^\a)}(h_n-\Psi^{\a,n}_f(\t))+\Psi^{\a,n}_f(t)\big]\f_n(x),\q0\leq t\leq\t.
\eeq
\et
\bpf
From Definition \ref{def-sol-backward-problem} and the solution representation \eqref{forward-sol-rep} of the solution of \eqref{forward-pde},  we observe that if $u(\cdot,t)$ as defined in \eqref{sol-rep-bacward-all-time}, belongs to $L^2(\O)$ for all $0\leq  t\leq  \t$, then it is indeed a solution of the backward problem \eqref{gov-pde}. Thus, it remains to ensure that $u(\cdot,t)\in L^2(\O)$ for $0\leq t\leq \t.$

Towards this, for $t\neq 0,$ by invoking Theorem \ref{prop_Ea_upper_bound_lower_bound}, we first observe that 
\beq\label{int-bound-interms-of-t}
\frac{\Ea(-\l_n^2t^\a)}{\Ea(-\l_n^2\t^\a)}\leq \frac{C_2(1+\l_1^2\t^\a)}{C_1\l_1^2}\frac{1}{t^\a},\eeq
 which is independent of $n.$  Hence by Propositions \ref{psi-alphan-difference-estimate} (i) and \ref{combined-exact-data} (i), it follows that $u(\cdot,t)\in L^2(\O).$
 
 For $t=0,$ observe that $\Ea(0)=1$ and $\Psi_{f}^{\a,n}(0)=0.$ Thus, $u(\cdot,0)\in L^2(\O)$ if and only if $$\sum_{n=1}^\infty \left(\frac{h_n-\Psi_{f}^{\a,n}(\t)}{\Ea(-\l_n^2\t^\a)}\right)^2<\infty.$$
\epf

\bt\label{stability-for-t-non-zero}
For $h,\tilde{h}\in L^2(\O),$ $f,\tilde{f}\in L^\infty(0,\t;L^2(\O)),$ let $u$ and $\tilde{u}$ be the corresponding solutions for the backward problem associated with \eqref{gov-pde}. Then for $0<t<\t,$ there exists a constant $C>0,$ depending on $C_1,C_2,\a,\t,\l_1,$ such that $$\|\tilde{u}(t)-u(t)\|^2\leq \frac{4C^2}{t^{2\a}}(\|h-h^\d\|^2+\theta\|f-f^\d\|^2_{L^\infty(0,\t;L^2(\O))})+2\theta\|f-f^\d\|^2_{L^\infty(0,\t;L^2(\O))}.$$
\et
\bpf
From \eqref{sol-rep-bacward-all-time}, we have
\beqarray
 \|u(t)-\tilde{u}(t)\|^2&\leq &2\sum_{n=1}^\infty |\frac{\Ea(-\l_n^2t^\a)}{\Ea(-\l_n^2\t^\a)}(h_n-\tilde{h}_n-\Psi_{f}^{\a,n}(\t)+\Psi_{\tilde{f}}^{\a,n}(\t))|^2+|\Psi_{f}^{\a,n}(t)-\Psi_{\tilde{f}}^{\a,n}(t)|^2\\
 &\underbrace{\leq }_{\eqref{int-bound-interms-of-t}}& 2 \left(\frac{C_2(1+\l_1^2\t^\a)}{C_1\l_1^2}\frac{1}{t^\a}\right)^2\sum_{n=1}^\infty |h_n-\tilde{h}_n-\Psi_{f}^{\a,n}(\t)+\Psi_{\tilde{f}}^{\a,n}(\t)|^2+2\sum_{n=1}^\infty|\Psi_{f}^{\a,n}(t)-\Psi_{\tilde{f}}^{\a,n}(t)|^2\\
 &=&2 \left(\frac{C_2(1+\l_1^2\t^\a)}{C_1\l_1^2}\frac{1}{t^\a}\right)^2\|\Upsilon_{h,f}-\Upsilon_{\tilde{h},\tilde{f}}\|^2+2\sum_{n=1}^\infty |\Psi_{f}^{\a,n}(t)-\Psi_{\tilde{f}}^{\a,n}(t)|^2,
\eeqarray
where recall that $\Upsilon_{h,f}$ is as in Proposition \ref{combined-exact-data}. Now the proof follows from Propositions \ref{psi-alphan-difference-estimate} and \ref{combined-exact-data}.
\epf
\brem\label{rem-stability-for-t-non-zero}
Theorem \ref{stability-for-t-non-zero} shows that, for $0<t<\t,$ the solution of the backward problem \eqref{gov-pde}, is stable with respect to the perturbation in the data $h$ and $f$.
\erem
The next question one would like to ask is, if $u$ is a solution of the backward problem \eqref{gov-pde}, then whether $u(0)$ is stable under perturbations in $h$ and $f.$ The answer is in the negative as is seen in the following example.
\beg{\bf (Ill-posed)}\label{eg-illposed}
For $n\in \N$ and $(x,t)\in \O\times (0,\t),$ let $f(x,t)=0=\tilde{f}^n(x,t),$ $h(x)=\f_1(x)$ and $\tilde{h}^n(x)=h(x)+\frac{1}{\l_n}\f_n(x).$ Corresponding to the source function $f$ and the final value $h,$ let $u(0)$ denote the initial value associated with the backward problem \eqref{gov-pde}, and for the source function $\tilde{f}$ and the final value $\tilde{h}^n,$ let $\tilde{u}^n(0)$ denote the initial value associated with the backward problem \eqref{gov-pde}. Then from the representation \eqref{sol-rep-bacward-all-time}, it follows that $u(0)=\frac{1}{\Ea(-\l_1^2\t^\a)}\f_1$ and $\tilde{u}^n(0)=\frac{1}{\Ea(-\l_1^2\t^\a)}\f_1+\frac{1}{\l_n\Ea(-\l_n^2\t^\a)}\f_n.$ Now we observe that $\|h-\tilde{h}^n\|=\frac{1}{\l_n}\to 0$ as $n\to \infty$ but $$\|u(0)-\tilde{u}^n(0)\|=\frac{1}{\l_n\Ea(-\l_n^2\t^\a)}\underbrace{\geq}_{Theorem \ref{prop_Ea_upper_bound_lower_bound}} \frac{\t^\a\Gamma(1-\a)}{C_2}\l_n\to \infty$$ as $n\to \infty.$
\eeg
We now obtain an equivalent formulation of the inverse problem of retrieving the initial value from the final value and the source function. 

Let $h\in L^2(\O)$, $f\in L^\infty(0,\t;L^2(\O))$ and $\Upsilon_{h,f}$ be as in Proposition \ref{combined-exact-data}. Let $u$ be the solution of the backward problem \eqref{gov-pde}. Then from \eqref{g-h-f-relation} and Proposition \ref{combined-exact-data} (i), we observe that the inverse problem of recovering the initial value $g=u(\cdot,0)$ from the knowledge of $h=u(\cdot,\t)$ and $f,$ is same as solving the linear operator equation 
\beq\label{op-eq}
Tg =\Upsilon_{h,f},
\eeq
where $T:L^2(\O)\to L^2(\O)$ is defined by
\beq\label{operator-def}
T\phi=\sum_{n=1}^\infty \Ea(-\l_n^2\t^\a)\<\phi,\f_n\>\f_n,\q\phi\in L^2(\O).
\eeq
From Theorem \ref{prop_Ea_upper_bound_lower_bound} it follows that $\Ea(-\l_n^2\t^\a)\sim\frac{1}{\l_n^2}$ and thus it follows that $T:L^2(\O)\to L^2(\O)$ is a compact linear operator of infinite rank. Therefore, solving \eqref{op-eq} is an ill-posed problem. Moreover, if $$\sum_{n=1}^\infty \left(\frac{\Upsilon_{h,f}^{\a,n}}{\Ea(-\l_n^2\t^\a)}\right)^2<\infty$$ then \eqref{op-eq}
 has a unique solution given by $$g=\sum_{n=1}^\infty g_n\f_n,$$ where $g_n$ is as defined in \eqref{g-h-f-relation} and $\Upsilon_{h,f}^{\a,n}$ is as defined in \eqref{psifan}.

We now recall definitions of some Sobolev spaces that are relevant in our analysis.
For $p\geq 0,$ let $$H_p=\big\{\phi\in L^2(\O):\,\sum_{n=1}^\infty \l_n^{2p}|\<\phi,\f_n\>|^2<\infty\big\}.$$ It is known that $H_p$ is a Hilbert space \w.r.t. the inner-product
$$\<\phi,\psi\>:=\sum_{n=1}^\infty \l_n^{2p}\<\phi,\f_n\>\<\f_n,\psi\>,\q\phi,\psi\in H_{p},$$ and the corresponding norm is given by $$\|\phi\|_{H_p}^2=\sum_{n=1}^\infty \l_n^{2p}|\<\phi,\f_n\>|^2,\q \phi\in H_p.$$
Note that $H_0=L^2(\O)$ and $H_{1/2}=H^1_0(\O).$

For some $\varrho>0,$ we consider the following set
\beq\label{source-set}
S_{\varrho,p}:=\big\{\phi\in H_p:\,\|\phi\|_{H_p}\leq \varrho\big\}.
\eeq
These sets will be the {\it source sets} that are of interest in our analysis.

We end this section by proving a conditional stability estimate.
\bt\label{cond-stability}
Let $h\in L^2(\O),f\in L^\infty(0,\t;L^2(\O))$ and $\Upsilon_{h,f}$ be as in Proposition \ref{combined-exact-data}. Let $g$ be the unique solution of \eqref{op-eq}. If $g\in S_{\varrho,p}$ for some $\varrho,p>0,$ then there exists a constant $C_3,$ depending on $\a,\t$ and $C_1,$ such that $$\|g\|\leq C_3^{p/(p+2)}\|g\|_{H_{p}}^{2/(p+2)}\|\Upsilon_{h,f}\|^{p/(p+2)}.$$
\et
\bpf
Since $g\in L^2(\O)$ is a solution of \eqref{op-eq}, we have $g=\sum_{n=1}^\infty g_n\f_n$, where $g_n$ are as in \eqref{g-h-f-relation}, that is, $$g_n=\frac{\Upsilon_{h,f}^{\a,n}}{\Ea(-\l_n^2\t^\a)}=\frac{h_n-\Psi^{\a,n}_{f}(\t)}{\Ea(-\l_n^2\t^\a)},\q n\in\N.$$
Thus,
\beqarray
\|g\|^2&=&\sum_{n=1}^\infty |\frac{\Upsilon_{h,f}^{\a,n}}{\Ea(-\l_n^2\t^\a)}|^2=\sum_{n=1}^\infty |\Upsilon_{h,f}^{\a,n}|^{2p/(p+2)}\,\left(\frac{|\Upsilon_{h,f}^{\a,n}|^2}{|\Ea(-\l_n^2\t^\a)|^{(p+2)}}\right)^{2/(p+2)}\\
&\leq & \left(\sum_{n=1}^\infty |\Upsilon_{h,f}^{\a,n}|^2\right)^{p/(p+2)}\left(\sum_{n=1}^\infty \frac{|\Upsilon_{h,f}^{\a,n}|^2}{|\Ea(-\l_n^2\t^\a)|^{(p+2)}}\right)^{2/(p+2)}.
\eeqarray
Now, from Theorem \ref{prop_Ea_upper_bound_lower_bound} we have $$\frac{1}{\Ea(-\l_n^2\t^\a)}\leq C_3\l_n^2,$$ where $C_3=\frac{(1+\l_1^2\t^\a)\Gamma(1-\a)}{C_1\l_1^2}.$ Thus,
using the fact that $$\frac{|\Upsilon_{h,f}^{\a,n}|^2}{|\Ea(-\l_n^2\t^\a)|^{(p+2)}}=\frac{g_n^2}{|\Ea(-\l_n^2\t^\a)|^{p}},$$ we have
$$\|g\|^2\leq C_3^{2p/(p+2)}\left(\sum_{n=1}^\infty |\Upsilon_{h,f}^{\a,n}|^2\right)^{p/(p+2)}\left(\sum_{n=1}^\infty \l_n^{2p}g_n^2\right)^{2/(p+2)},$$
that is,
$$\|g\|^2\leq C_3^{2p/(p+2)}\|\Upsilon_{h,f}\|^{2p/(p+2)}\|g\|_{H_{p}}^{4/(p+2)}.$$
Now the proof follows.
\epf
\section{Quasi-boundary value method and its modification}\label{sec-quasi-bound}
We have seen in the preceding section that the considered inverse problem of retrieving the initial value from the knowledge of final value and the source term, is an ill-posed problem in the sense that a small perturbations in the data $h$ and $f$ may lead to a large deviation in the sought solution. But, in practice, we shall have knowledge of the data that are corrupted by some noise. Thus, we have to employ some regularization method to obtain stable approximations. In this section we shall employ the quasi-boundary value method and its modifications, and obtain error estimates for apriori as well as aposteriori parameter choice strategies.

For $\b>0$ and $q=\N\cup\{0\},$ we consider the following equation
\beq\label{quasi-bound-pde}
\begin{cases}
\partial^\a_t v+\Delta^2 v=f,&\q\text{in}\q\O\times (0,\t),\\
v=0=\Delta v,&\q\text{on}\q\partial\O\times (0,\t),\\
v(\t)+\b (-\Delta)^q v(0)=h,&\q\text{in}\q\O,
\end{cases}
\eeq
where $h\in L^2(\O)$ and $f\in L^\infty(0,\t;L^2(\O)).$ Note that $\b$ will play the role of regularization parameter.

By separation of variable method, we know that the solution of \eqref{quasi-bound-pde} is of the form
$$u_\b(x,t)=\sum_{n=1}^\infty \big[A_n \Ea(-\l_n^2t^\a)+\int_0^t (t-s)^{\a-1}\Eaa(-\l_n^2(t-s)^\a)f_n(s)\,ds\big]\f_n(x),$$ where $A_n$'s are constants to be determined. Now using the condition $u_\b(\t)+\b(-\Delta)^qu_\b(0)=h,$ we obtain
$$A_n=\frac{h_n-\Psi_{f}^{\a,n}(\t)}{\Ea(-\l_n^2\t^\a)+\b\l_n^q},\q n\in\N,$$ where $\Psi_{f}^{\a,n}(\t)$ is as defined in \eqref{psifan}. Therefore, the solution $u_\b$ of \eqref{quasi-bound-pde} is given by
\beq\label{sol-rep-quasi-bound-pde}
u_\b(x,t)=\sum_{n=1}^\infty \big[\frac{(h_n-\Psi_{f}^{\a,n}(\t))\Ea(-\l_n^2t^\a)}{\Ea(-\l_n^2\t^\a)+\b\l_n^q}+\Psi_{f}^{\a,n}(t)\big]\f_n(x)
\eeq
Now taking the noisy data $h^\d$ and $f^\d$ in place of the exact data $h$ and $f$, respectively, the solution $u_\b^\d$ of \eqref{quasi-bound-pde} is given by
\beq\label{sol-rep-quasi-bound-pde-noisy}
u_\b^\d(x,t)=\sum_{n=1}^\infty \big[\frac{(h^\d_n-\Psi_{f^\d}^{\a,n}(\t))\Ea(-\l_n^2t^\a)}{\Ea(-\l_n^2\t^\a)+\b\l_n^q}+\Psi_{f^\d}^{\a,n}(t)\big]\f_n(x),
\eeq
where $\Psi_{f^\d}^{\a,n}(\t)$ is as defined in \eqref{psifan}.
\bl\label{sup-abeta}
Let $f\in L^\infty(0,\t;L^2(\O))$ and $h\in L^2(\O)$. For $\b>0$, let $A_\b(n):=\frac{1}{\Ea(-\l_n^2\t^\a)+\b\l_n^q},\,\,n\in\N$ and $q\in\N\cup\{0\}.$ The following holds:
\ben
\i[(i)] There exists a constant $C_4>0,$ depending on $\a,\t,\l_1,C_1,$ such that $$\sup_n A_\b(n)\leq \begin{cases}
\frac{1}{\b},&\q q=0,\\
(2C_4/q)^{2/(q+2)}\frac{q}{C_4(1+q)}\frac{1}{\b^{2/(q+2)}},&\q q\neq 0.
\end{cases}$$
\i[(ii)] For $q\in\N\cup\{0\},$ we have $$\sum_{n=1}^\infty \big|\frac{h_n-\Psi_{f}^{\a,n}(\t)}{\Ea(-\l_n^2\t^\a)+\b\l_n^q}\big|^2\leq 2(\sup_n A_\b(n))^2(\|h\|^2+\theta\|f\|^2_{L^\infty(0,\t;L^2(\O))})<\infty,$$ where $\theta$ is as in Proposition \ref{combined-exact-data}.
\een
\el
\bpf
{\bf (i).} Let $q=0.$ Since $\Ea(-\l_n^2\t^\a)>0,$ it follows that $A_\b(n)\leq \frac{1}{\b},$ for all $n.$

Let $q\neq 0.$ From Theorem \ref{prop_Ea_upper_bound_lower_bound}, we have
$$\Ea(-\l_n^2\t^\a)\geq \frac{C_4}{\l_n^2},$$ where $C_4:= \frac{C_1}{\Gamma(1-\a)}\frac{\l_1^2}{1+\l_1^2\t^\a}.$ Thus, we have
$$A_\b(n)\leq \frac{\l_n^2}{C_4+\b \l_n^{q+2}},\q n\in\N.$$ Now the result follows from Lemma \ref{calculus} (i).

\noi
{\bf (ii).}
The proof follows from Proposition \ref{combined-exact-data} and Lemma \ref{sup-abeta} (i).
\epf
\brem\label{existence-initial-quasi-boundary}
Let $h\in L^2(\O)$ and $f\in L^\infty(0,\t;L^2(\O)).$ It is worth to note that for every $\b>0$ and $q\in \N\cup\{0\},$ we have $u_{\b}(t)\in L^2(\O)$ for all $0\leq t\leq \t,$ where $u_\b(t)$ is as defined in \eqref{sol-rep-quasi-bound-pde}. Indeed, for $0<t\leq \t,$ the assertion follows by the arguments employed in the proof of Theorem \ref{nec-suff-sond-backward-problem}. For $t=0,$ the assertion follows by Lemma \ref{sup-abeta} (ii).
\erem
Thus, it follows that for every $\b>0,q\in\N\cup\{0\}, h\in L^2(\O)$ and $f\in L^\infty(0,\t;L^2(\O)),$ the equation \eqref{quasi-bound-pde} has a unique solution $u_\b$ given by \eqref{sol-rep-quasi-bound-pde}. In particular, we have 
\beq\label{initial-value-quasi-bound-exact}
u_\b(0)=\sum_{n=1}^\infty \big[\frac{h_n-\Psi_f^{\a,n}(\t)}{\Ea(-\l_n^2\t^\a)+\b\l_n^q}\big]\f_n\in L^2(\O).
\eeq
Note that these $u_\b(0)$ are the candidates for the regularized solutions in this section.

For $h^\d\in L^2(\O)$ and $f^\d\in L^\infty(0,\t;L^2(\O)),$ from \eqref{sol-rep-quasi-bound-pde-noisy}, we have
\beq\label{initial-value-quasi-bound-noisy}
u^\d_\b(0)=\sum_{n=1}^\infty \big[\frac{h^\d_n-\Psi_{f^\d}^{\a,n}(\t)}{\Ea(-\l_n^2\t^\a)+\b\l_n^q}\big]\f_n\in L^2(\O)
\eeq

\bt\label{stability-initial-quasi-bound}
Let $\theta$ be as in Proposition \ref{combined-exact-data}, $C_4$ be as in Lemma \ref{sup-abeta} and $q\in \N\cup\{0\}$. For $\d>0,$ let $f,f^\d\in L^\infty(0,\t;L^2(\O))$ and $h,h^\d\in L^2(\O)$ be such that $$\|h-h^\d\|^2+\theta\|f-f^\d\|^2_{L^\infty(0,\t:L^2(\O))}\leq \d^2.$$ For $\b>0,$ let $u_\b(0)$ and $u_\b^\d(0)$ be as defined in \eqref{initial-value-quasi-bound-exact} and \eqref{initial-value-quasi-bound-noisy}, respectively. Then there exists a constant $C_5>0,$ depending on $C_4,q,$ such that $$\|u_\b^\d(0)-u_\b(0)\|\leq 
\begin{cases}
\sqrt{2}\frac{\d}{\b},&\q q=0,\\
C_5\frac{\d}{\b^{2/(q+2)}},&\q q\neq 0.
\end{cases}
$$
\et
\bpf
We have 
\beqarray
\|u_\b^\d(0)-u_\b(0)\|^2&=&\sum_{n=1}^\infty \left|\frac{h_n^\d-h_n-\Psi_{f^\d}^{\a,n}(\t)+\Psi_{f}^{\a,n}(\t)}{\Ea(-\l_n^2\t^\a)+\b\l_n^q}\right|^2\\
&\leq & \left(\sup_n \frac{1}{\Ea(-\l_n^2\t^\a)+\b\l_n^q}\right)^2\sum_{n=1}^\infty |h_n^\d-h_n-\Psi_{f^\d}^{\a,n}(\t)+\Psi_{f}^{\a,n}(\t)|^2\\
&\leq &2\left(\sup_n \frac{1}{\Ea(-\l_n^2\t^\a)+\b\l_n^q}\right)^2\sum_{n=1}^\infty (|h_n-h_n^\d|^2+|\Psi_{f^\d}^{\a,n}-\Psi_f^{\a,n}|^2)\\
&\underbrace{\leq}_{Lemma \ref{psi-alphan-difference-estimate}} &2\left(\sup_n \frac{1}{\Ea(-\l_n^2\t^\a)+\b\l_n^q}\right)^2\sum_{n=1}^\infty (|h_n-h_n^\d|^2+ \frac{\|f-f^\d\|^2_{L^\infty(0,\t;L^2(\O))}}{\l_n^4})\\
&=&2\left(\sup_n \frac{1}{\Ea(-\l_n^2\t^\a)+\b\l_n^q}\right)^2 (\|h-h^\d\|^2+\theta\|f-f^\d\|^2_{L^\infty(0,\t;L^2(\O))})\\
&\leq & 2\left(\sup_n \frac{1}{\Ea(-\l_n^2\t^\a)+\b\l_n^q}\right)^2 \d^2.
\eeqarray
Now the proof follows from Lemma \ref{sup-abeta} (i) for $q=0$, and by taking $C_5=\sqrt{2}(2C_4/q)^{2/(q+2)}\frac{q}{C_4(1+q)}$ for $q\neq 0.$
\epf
\brem
From Remark \ref{existence-initial-quasi-boundary} and Theorem \ref{stability-initial-quasi-bound}, it follows that retrieving $u_\b(0)$ from the knowledge of $f,h$ associated with \eqref{quasi-bound-pde}, is an well-posed problem.
\erem
For $\varrho,p>0,$ recall the definition of the source set $S_{\varrho,p}$ as given in \eqref{source-set}.
\bt\label{estimate-exact-regularized}
Let $q\in \N\cup\{0\}$ and for $\b>0,$ let $u_\b(0)$ be as in \eqref{initial-value-quasi-bound-exact}, $C_4$ be as in Lemma \ref{sup-abeta} and $g$ be the unique solution of \eqref{op-eq}. If $g\in S_{\varrho,p},$ for some $\varrho,p>0,$ then there exist constants $C_6>0,$ depending on $C_4, q, p,$ and $C_7>0,$ depending on $C_4,p,q,\l_1,$ such that $$\|u_\b(0)-g\|\leq 
\begin{cases}
C_6 \varrho \b^{p/(q+2)},&\q 0<p<q+2,\\
C_7\varrho\b,&\q p\geq q+2.
\end{cases}$$
\et
\bpf
Recall that $g=\sum_{n=1}^\infty g_n\f_n,$ where $g_n$ is as in \eqref{g-h-f-relation}. Thus, we have
\beqarray
\|u_\b(0)-g\|^2&=& \sum_{n=1}^\infty \left|\left(\frac{1}{\Ea(-\l_n^2\t^\a)+\b\l_n^q}-\frac{1}{\Ea(-\l_n^2\t^\a)}\right)(h_n-\Psi_f^{\a,n}(\t))\right|^2\\
&=& \sum_{n=1}^\infty \left(\frac{\b\l_n^{q-p}}{\Ea(-\l_n^2\t^\a)+\b\l_n^q}\right)^2\left(\frac{\l_n^{p}(h_n-\Psi_f^{\a,n}(\t))}{\Ea(-\l_n^2\t^\a)}\right)^2\\
&\leq& (\sup_n B_\b(n))^2\sum_{n=1}^\infty \l_n^{2p}g_n^2\\
&\leq & (\sup_n B_\b(n))^2\|g\|^2_{H_{p}},
\eeqarray
where $B_\b(n)=\frac{\b\l_n^{q-p}}{\Ea(-\l_n^2\t^\a)+\b\l_n^q},\,\,n\in\N.$

Now from Theorem \ref{prop_Ea_upper_bound_lower_bound}, we observe that $$\Ea(-\l_n^2\t^\a)+\b\l_n^q\geq \frac{C_4+\b\l_n^{q+2}}{\l_1^2},$$ where $C_4=\frac{C_1\l_1^2}{(1+\l_1^2\t^\a)\Gamma(1-\a)}.$

\noi
{\bf Case 1.} Let $0<p<q+2.$

Then by Lemma \ref{calculus} (ii), it follows that $$B_\b(n)\leq C_6 \b^{p/(q+2)},$$ where $C_6=\left(\frac{C_4(q+2-p)}{p}\right)^{(q+2-p)/(q+2)}\frac{p}{C_4(q+2)}.$

\noi
{\bf Case 2.} Let $p\geq q+2.$ 

Then $$B_\b(n)\leq \frac{\b\l_n^{q+2-p}}{C_4+\b\l_n^{q+2}}\leq \frac{\b}{C_4\l_n^{p-q-2}}\leq C_7 \b,$$where $C_7=\frac{1}{C_4\l_1^{p-q-2}}.$

Now the proof follows.
\epf
Now we are in a position to state one of the main result of this section whose proof follows by combining the estimates obtained in Theorems \ref{stability-initial-quasi-bound} and \ref{estimate-exact-regularized}, respectively. This result provides the rates of convergence for the approximations by regularized solutions obtained by modified quasi-boundary value method for all $q\in \N\cup\{0\}$ under an apriori parameter choice strategy.
\bt\label{quasi-bound-apriori-rate}
For $\b,\d>0$ and $q\in \N\cup\{0\},$ let $h,h^\d,f,f^\d,u_\b(0),u_\b^\d(0)$ be as in Theorem \ref{stability-initial-quasi-bound}. Let $p,\varrho, g$ be as in Theorem \ref{estimate-exact-regularized}. Let $C_5,C_6,C_7$ be as in Theorems \ref{stability-initial-quasi-bound} and \ref{estimate-exact-regularized}, respectively. Then we have the following:
\ben
\i[(i)] For $q=0,\, 0<p<2,$ and the choice $\b\sim (\d/\varrho)^{2/(p+2)},$ there exists a constant $\tilde{C}_8>0,$ depending on $C_6,p,$ such that $$\|u_\b^\d(0)-g\|\leq \tilde{C}_8 \varrho^{2/(p+2)}\d^{p/(p+2)}.$$
\i[(ii)] For $q=0,\,p\geq 2,$ and the choice $\b\sim (\d/\varrho)^{1/2},$ there exists a constant $\tilde{C}_9,$ depending on $C_7,$ such that $$\|u_\b^\d(0)-g\|\leq \tilde{C}_9 \varrho^{1/2}\d^{1/2}.$$
\i[(iii)] For $q\neq 0,\,0<p<q+2,$ and the choice $\b\sim (\d/\varrho)^{(q+2)/(p+2)},$ there exists a constant $C_8>0,$ depending on $C_5,C_6,p,q$, such that $$\|u_\b^\d(0)-g\|\leq C_8 \varrho^{2/(p+2)}\d^{p/(p+2)}.$$
\i[(iv)] For $q\neq 0,\,p\geq q+2$ and the choice $\b\sim (\d/\varrho)^{(q+2)/(q+4)},$ there exists a constant $C_9>0,$ depending on $C_5,C_7,q$, such that $$\|u_\b^\d(0)-g\|\leq C_9 \varrho^{2/(q+4)}\d^{(q+2)/(q+4)}.$$
\een
\et
In Theorem \ref{quasi-bound-apriori-rate} we have obtained the rates by choosing the regularization parameter $\b$ in a certain manner that requires the apriori knowledge of $p.$ This means that one has to apriori know the smoothness of the unknown $g.$ But that is not always feasible and is too restrictive to assume. In order to overcome this, we provide an aposteriori parameter choice strategy in the subsequent analysis.  
\subsection{Error estimates with aposteriori parameter choice}
For $h\in L^2(\O)$ and $f\in L^2(\O),$ let $\Upsilon_{h,f}$ be as in Proposition \ref{combined-exact-data}, that is, $\Upsilon_{h,f}=\sum_{n=1}^\infty \Upsilon_{h,f}^{\a,n}\f_n,$ where $\Upsilon_{h,f}^{\a,n}$ is as defined in \eqref{psifan}. Also, recall that the considered inverse problem is formulated as solving an operator equation $Tg=\Upsilon_{h,f},$ where $T$ is as defined in \eqref{operator-def}.  

We now state a result which has been used frequently in the literature. Since the context in which we are using is new, we include only the key steps in the proof for the sake of completeness.
\bl\label{apost-discr-function-quasi-bound}
For $\b,\d>0$ and $q\in \N\cup\{0\},$ let $u_\b^\d(0)$ be as in \eqref{initial-value-quasi-bound-noisy} and $\Upsilon_{h^\d,f^\d}$ be as defined above with $\|\Upsilon_{h^\d,f^\d}\|\neq 0.$ Let $$\Phi(\b)=\|Tu_\b^\d(0)-\Upsilon_{h^\d,f^\d}\|,\q\b>0.$$ Then the following holds:
\ben
\i[(i)] $\Phi$ is continuous and increasing.
\i[(ii)] $\lim_{\b\to 0}\Phi(\b)=0$ and $\lim_{\b\to \infty}\Phi(\b)=\|\Upsilon_{h^\d,f^\d}\|.$
\een
\el
\bpf
{\bf (i).} The continuity follows easily by observing that $$(\Phi(\b))^2=\sum_{n=1}^\infty \frac{\b^2\l_n^{2q}(\Upsilon_{h^\d,f^\d}^{\a,n})^2}{(\Ea(-\l_n^2\t^\a)+\b\l_n^q)^2}.$$ In order to show that $\Phi$ is increasing, we first observe that for $b,c>0,$ the function $\frac{bt}{bt+c}$ is increasing for $t\geq 0.$ Since $\|\Upsilon_{h^\d,f^\d}\|\neq 0,$ there exists some $m\in\N,$ such that $\Upsilon_{h^\d,f^\d}^{\a,m}\neq 0.$ Using these facts, it follows that $\Phi$ is indeed increasing.

\noi
{\bf (ii).} From the expression of $(\Phi(\b))^2$ it follows that $\lim_{\b\to 0}\Phi(\b)=0.$ Next, from Theorem \ref{prop_Ea_upper_bound_lower_bound}, we obtain $$\frac{\b^2\l_n^{2q}}{(\Ea(-\l_n^2\t^\a)+\b\l_n^q)^2}\geq \frac{1}{1+\tilde{C}/(\b\l_1^q(1+\l_1^2\t^\a))},$$ where $\tilde{C}=\frac{C_2}{\Gamma(1-\a)}.$ This implies that $$(\Phi(\b))^2\geq \frac{1}{1+\tilde{C}/(\b\l_1^q(1+\l_1^2\t^\a))}\|\Upsilon_{h^\d,f^\d}\|^2.$$ However, we also observe that $$(\Phi(\b))^2\leq \|\Upsilon_{h^\d,f^\d}\|^2.$$ Now the proof follows by taking the limit as $\b\to \infty.$
\epf
\bl\label{reg-parameter-bound-aposteriori-quasi-bound} 
Let $g$ be the solution of \eqref{op-eq} and suppose that $g\in S_{\varrho,p}$ for some $\varrho,p>0.$ For $\d>0,$ let $\Phi$ be as defined in Lemma \ref{apost-discr-function-quasi-bound}. Let $\xi>\sqrt{2}$ be given and $\nu\in(0,1)$ be fixed. Then for every $\d>0$ small enough, with $(\xi+\sqrt{2})\d^\nu\leq \|\Upsilon_{h,f}\|$ there exist unique $\b_\d>0$ and $\tilde{\b}_\d>0$ such that $$\Phi(\b_\d)=\xi\d,\q\Phi(\tilde{\b}_\d)=\xi\d^\nu,$$ and the following result holds:
\ben
\i[(i)] For $q\neq 0$ and $0<p<q,$ there exists a constant $C_{11}>0,$ depending on $C_1,C_2,\a,\t,\l_1,p,q,$ such that $$\frac{1}{\b_\d}\leq \left(\frac{C_{11}}{\xi-\sqrt{2}}\right)^{(q+2)/(p+2)}\left(\frac{\varrho}{\d}\right)^{(q+2)/(p+2)}.$$
\i[(ii)] For $q\neq 0$ and $p\geq q,$ there exists a constant $C_{12}>0,$ depending on $C_1,C_2,\a,\t,\l_1,p,q,$ such that $$\frac{1}{\b_\d}\leq \frac{C_{12}}{\xi-\sqrt{2}}\frac{\varrho}{\d}.$$
\i[(iii)] For $q=0$ we have $$\frac{1}{\tilde{\b}_\d}\leq \frac{1}{\l_1^{2p}(\xi-\sqrt{2})}\frac{\varrho}{\d^\nu}.$$
\een
\el
\bpf
 From Proposition \ref{combined-exact-data} (ii), we first note that for $0<\d<1,$ 
 $$\|\Upsilon_{h,f}-\Upsilon_{h^\d,f^\d}\|\leq \sqrt{2}\d<\sqrt{2}\d^\nu.$$ Let $\xi>\sqrt{2}.$ Now the existence of unique $\b_\d,\tilde{\b}_\d>0$ satisfying $\Phi(\b_\d)=\xi\d$ and $\Phi(\tilde{\b}_\d)=\xi\d^\nu$ follows from Lemma \ref{apost-discr-function-quasi-bound} by observing that $\|\Upsilon_{h^\d,f^\d}\|\geq \|\Upsilon_{h,f}\|-\|\Upsilon_{h,f}-\Upsilon_{h^\d,f^\d}\|\geq \xi\d^\nu>\xi\d.$  Thus, we have
\beqarray
\xi\d &=& \Phi(\b_\d)= \|\sum_{n=1}^\infty \frac{\Ea(-\l_n^2\t^\a)\Upsilon_{h^\d,f^\d}^{\a,n}}{\Ea(-\l_n^2\t^\a)+\b_\d\l_n^q}\f_n-\sum_{n=1}^\infty \Upsilon_{h^\d,f^\d}^{\a,n}\f_n\|\\
&=& \|\sum_{n=1}^\infty \frac{\b_\d\l_n^q\Upsilon_{h^\d,f^\d}^{\a,n}}{\Ea(-\l_n^2\t^\a)+\b_\d\l_n^q}\f_n\|\\
&\leq & \|\sum_{n=1}^\infty \frac{\b_\d\l_n^q \Upsilon_{h,f}^{\a,n}}{\Ea(-\l_n^2\t^\a)+\b_\d\l_n^q}\f_n\|+\|\sum_{n=1}^\infty \frac{\b_\d\l_n^q(\Upsilon_{h^\d,f^\d}^{\a,n}-\Upsilon_{h,f}^{\a,n})}{\Ea(-\l_n^2\t^\a)+\b_\d\l_n^q}\f_n\|\\
&\leq & \|\sum_{n=1}^\infty \frac{\b_\d\l_n^q \Upsilon_{h,f}^{\a,n}}{\Ea(-\l_n^2\t^\a)+\b_\d\l_n^q}\f_n\| +\sqrt{2}\d.
\eeqarray
Thus, 
\beqarray
(\xi-\sqrt{2})^2\d^2 &\leq & \sum_{n=1}^\infty \left(\frac{\b_\d\l_n^q \Upsilon_{h,f}^{\a,n}}{\Ea(-\l_n^2\t^\a)+\b_\d\l_n^q}\right)^2 \\
&=& \sum_{n=1}^\infty (D_{\b_\d}(n))^2 \frac{\l_n^{2p}}{(\Ea(-\l_n^2\t^\a))^2}(\Upsilon_{h,f}^{\a,n})^2,
\eeqarray
where $D_{\b_\d}(n)=\frac{\b_\d\l_n^q\Ea(-\l_n^2\t^\a)}{\l_n^{p}(\Ea(-\l_n^2\t^\a)+\b_\d\l_n^q)}.$
Now, from Theorem \ref{prop_Ea_upper_bound_lower_bound}, it can be verified that
$$D_{\b_\d}(n)\leq C_{10}\frac{\b_\d\l_n^{q-p}}{C_4+\b_\d\l_n^{q+2}},$$ where $C_4=\frac{C_1}{\Gamma(1-\a)}\frac{\l_1^2}{1+\l_1^2\t^\a}$ and $C_{10}=\frac{C_2}{\Gamma(1-\a)\t^\a}.$

\noi
{\bf (i).} Let $ q\neq 0$ and $0<p<q.$

Then we have
\beqarray
C_4+\b_\d\l_n^{q+2}&\geq & \frac{q-p}{q+2}\b_\d\l_n^{q+2} + \frac{p+2}{q+2} C_4\\
&=& \frac{q-p}{q+2}\big[(\b_\d\l_n^{q+2})^{(q-p)/(q+2)}\big]^{(q+2)/(q-p)} + \frac{p+2}{q+2}\big[C_4^{(p+2)/(q+2)}\big]^{(q+2)/(p+2)}\\
&\geq & (\b_\d\l_n^{q+2})^{(q-p)/(q+2)}\,(C_4)^{(p+2)/(q+2)}\\
&=& C_4^{(p+2)/(q+2)}{\b_\d}^{(q-p)/(q+2)}\l_n^{q-p}.
\eeqarray
Thus, in this case, 
$$D_{\b_\d}(n)\leq C_{10}\frac{\b_\d}{C_4^{(p+2)/(q+2)}{\b_\d}^{(q-p)/(q+2)}}=C_{11}\,{\b_\d}^{(p+2)/(q+2)},$$ where $C_{11}=\frac{C_{10}}{C_4^{(p+2)/(q+2)}}.$ Therefore,
\beqarray
(\xi-\sqrt{2})^2\d^2 &\leq & C_{11}^2 {\b_\d}^{2(p+2)/(q+2)}\sum_{n=1}^\infty \frac{\l_n^{2p}}{(\Ea(-\l_n^2\t^\a))^2}(\Upsilon_{h,f}^{\a,n})^2\\
&=& C_{11}^2 {\b_\d}^{2(p+2)/(q+2)} \sum_{n=1}^\infty \l_n^{2p}g_n^2\\
&\leq &C_{11}^2 \varrho^2 {\b_\d}^{2(p+2)/(q+2)}.
\eeqarray
Thus, we have
$$\frac{1}{\b_\d}\leq \left(\frac{C_{11}}{\xi-\sqrt{2}}\right)^{(q+2)/(p+2)}\left(\frac{\varrho}{\d}\right)^{(q+2)/(p+2)}.$$

\noi
{\bf (ii).} Let $q\neq 0$ and $p\geq q.$

In this case, we have
$$D_{\b_\d}(n)\leq C_{10}\frac{\b_\d\l_n^{q-p}}{C_4+\b_\d\l_n^{q+2}}\leq \frac{C_{10}}{C_4}\frac{\b_\d}{\l_1^{p-q}}=C_{12}\b,$$ where $C_{12}=\frac{C_{10}}{C_4 \l_1^{p-q}}.$ Thus,
$$(\xi-\sqrt{2})^2\d^2\leq C_{12}^2 {\b_\d}^2\varrho^2,$$ and hence
$$\frac{1}{\b_\d}\leq \frac{C_{12}}{\xi-\sqrt{2}}\frac{\varrho}{\d}.$$

\noi
{\bf (iii).} Let $q=0.$ Recall that there exists a unique $\tilde{\b}_\d>0$ such that $\Phi(\tilde{\b}_\d)=\xi\d^\nu.$ Thus, we have
\beqarray
\xi \d^\nu &=& \|\sum_{n=1}^\infty \left(\frac{\Ea(-\l_n^2\t^\a)}{\Ea(-\l_n^2\t^\a)+\tilde{\b}_\d}\Upsilon_{h^\d,f^\d}^{\a,n}-\Upsilon_{h^\d,f^\d}^{\a,n}\right)\f_n\|=\|\sum_{n=1}^\infty \frac{\tilde{\b}_\d\Upsilon_{h^\d,f^\d}^{\a,n}}{\Ea(-\l_n^2\t^\a)+\tilde{\b}_\d}\f_n\|\\
&\leq &\|\sum_{n=1}^\infty \frac{\tilde{\b}_{\d}\Upsilon_{h,f}^{\a,n}}{\Ea(-\l_n^2\t^\a)+\tilde{\b}_{\d}}\f_n\|+\sqrt{2}\d.
\eeqarray
Therefore,
\beqarray
(\xi-\sqrt{2})^2\d^{2\nu}&\leq & \sum_{n=1}^\infty \frac{{\tilde{\b}_{\d}}^2(\Upsilon_{h,f}^{\a,n})^2}{(\Ea(-\l_n^2\t^\a)+\tilde{\b}_{\d})^2}\\
&=& \sum_{n=1}^\infty \frac{{\tilde{\b}_{\d}}^2(\Ea(-\l_n^2\t^\a))^2}{(\Ea(-\l_n^2\t^\a)+\tilde{\b}_{\d})^2} g_n^2\\
&=& \sum_{n=1}^\infty  \frac{{\tilde{\b}_{\d}}^2(\Ea(-\l_n^2\t^\a))^2}{\l_n^{2p}(\Ea(-\l_n^2\t^\a)+\tilde{\b}_{\d})^2}\l_n^{2p} g_n^2\\
&\leq & \frac{{\tilde{\b}_{\d}}^2}{\l_1^{2p}}\varrho^2,
\eeqarray
and hence
$$(\xi-\sqrt{2})\d^\nu\leq \frac{\varrho\tilde{\b}_{\d}}{\l_1^{p}}.$$ Therefore,
$$\frac{1}{\tilde{\b}_\d}\leq \frac{1}{\l_1^{p}(\xi-\sqrt{2})}\frac{\varrho}{\d^\nu}.$$
\epf
Thus, from Theorem \ref{stability-initial-quasi-bound} and Lemma \ref{reg-parameter-bound-aposteriori-quasi-bound}, we have the following result.
\bt\label{estimate-reg-noisy-aposteriori-quasi-bound}
Let $\xi, \nu \in (0,1),\d,\Phi,\b_\d$ and $\tilde{\b}_\d$ be as in Lemma \ref{reg-parameter-bound-aposteriori-quasi-bound}. Let $u^\d_{\b_\d}(0)$ and $u_{\b}^\d(0)$ be as in \eqref{initial-value-quasi-bound-noisy} and \eqref{initial-value-quasi-bound-exact}, respectively. Let $g$ be the solution of \eqref{op-eq} and $C_5, C_{11}, C_{12}$ be as in Theorem \ref{stability-initial-quasi-bound} and Lemma \ref{reg-parameter-bound-aposteriori-quasi-bound}. If $g\in S_{\varrho,p}$ for some $\varrho,p>0,$ then the following result holds:
\ben
\i[(i)] For $q\neq 0,\,0<p<q,$ there exists a constant $C_{14}>0,$ depending on $C_5,C_{11}, \xi,p,$ such that $$\|u_{\b_\d}^\d(0)-u_{\b_\d}\|\leq C_{14}\varrho^{2/(p+2)}\d^{p/(p+2)}.$$
\i[(ii)] For $q\neq 0,\,p\geq q,$ there exists a constant $C_{15}>0,$ depending on $C_5, C_{12}, \xi,q,$ such that $$\|u_{\b_\d}^\d(0)-u_{\b_\d}(0)\|\leq C_{15}\varrho^{2/(q+2)}\d^{q/(q+2)}.$$
\i[(iii)] For $q=0,$ there exists a constant $C_{16}>0,$ depending on $C_5,\xi,\l_1,p,$ such that $$\|u_{\tilde{\b}_\d}^\d(0)-u_{\tilde{\b}_\d}(0)\|\leq C_{16}\varrho \d^{1-\nu}.$$
\een
\et
We now obtain estimate for $\|u_{\b_\d}(0)-g\|$ and $\|u_{\tilde{\b}_\d}(0)-g\|.$
\bt\label{estimate-reg-exact-aposteriori-quasi-bound}
Let $\d,\nu,\xi,\b_\d,\tilde{\b}_\d$ be as in Theorem \ref{estimate-reg-noisy-aposteriori-quasi-bound}, $u_{\b_\d}(0)$ be as in \eqref{initial-value-quasi-bound-exact} and $C_3$ be as in Theorem \ref{cond-stability}. Let $g$ be the solution of \eqref{op-eq}. If $g\in S_{\varrho,p},$ for some $\varrho,p>0$ then there exists a constant $C'_{13}>0,$ depending on $C_3,p,\xi,$ such that the following holds:
\ben
\i[(i)] For $q\neq 0,$  $$\|g-u_{\b_\d}(0)\|\leq C'_{13} \varrho^{2/(p+2)}\d^{p/(p+2)}.$$
\i[(ii)] For $q=0,$  $$\|g-u_{\tilde{\b}_\d}(0)\|\leq C'_{13} \varrho^{2/(p+2)}\d^{p\nu/(p+2)}.$$
\een
\et
\bpf From the definition of $T$ (see \eqref{operator-def}), we have
$$Tu_{\b_\d}(0)-Tu^\d_{\b_\d}(0)=\sum_{n=1}^\infty \frac{\Ea(-\l_n^2\t^\a)}{\Ea(-\l_n^2\t^\a)+\b_\d\l_n^q}(\Upsilon_{h^\d,f^\d}^{\a,n}-\Upsilon_{h,f}^{\a,n})\f_n,$$
and hence
$$\|T u_{\b_\d}(0)-Tu^\d_{\b_\d}(0)\|^2=\sum_{n=1}^\infty \big|\frac{\Ea(-\l_n^2\t^\a)}{\Ea(-\l_n^2\t^\a)+\b_\d\l_n^q}\big|^2|\Upsilon_{h^\d,f^\d}^{\a,n}-\Upsilon_{h,f}^{\a,n}|^2\leq 2\d^2.$$

\noi
{\bf (i).} Let $q\neq 0.$

In this case, we have
\beqarray
\|Tg-Tu_{\b_\d}(0)\|&=&\|\Upsilon_{h,f}-Tu_{\b_\d}(0)\|\\
&\leq &\|Tu_{\b_\d}(0)-Tu^\d_{\b_\d}(0)\|+\|Tu^\d_{\b_\d}(0)-\Upsilon_{h^\d,f^\d}\|+\|\Upsilon_{h,f}-\Upsilon_{h^\d,f^\d}\|\\
&\leq & \sqrt{2}\d+\xi\d +\sqrt{2}\d \\
&=& (\xi+2\sqrt{2})\d.
\eeqarray
Now, recall that 
\beqarray
g-u_{\b_\d}(0)&=& \sum_{n=1}^\infty \frac{\Upsilon_{h,f}^{\a,n}}{\Ea(-\l_n^2\t^\a)}\f_n-\sum_{n=1}^\infty \frac{\Upsilon_{h,f}^{\a,n}}{\Ea(-\l_n^2\t^\a)+\b\l_n^q}\f_n\\
&=&\sum_{n=1}^\infty \frac{\b\l_n^qg_n}{\Ea(-\l_n^2\t^\a)+\b\l_n^q}\f_n.
\eeqarray
Thus,
$$\sum_{n=1}^\infty \l_n^{2p}|\<g-u_{\b_\d}(0),\f_n\>|^2=\sum_{n=1}^\infty \l_n^{2p}\left(\frac{\b\l_n^qg_n}{\Ea(-\l_n^2\t^\a)+\b\l_n^q}\right)^2\leq \sum_{n=1}^\infty \l_n^{2p}g_n^2\leq \varrho^2.$$
This shows that $g-u_{\b_\d}(0)\in H_{p}$ and $\|g-u_{\b_\d}(0)\|_{H_{p}}\leq \varrho.$ Thus, from Theorem \ref{cond-stability}, we have
$$\|g-u_{\b_\d}(0)\|\leq C_3^{p/(p+2)}\varrho^{2/(p+2)}(\xi+2\sqrt{2})^{p/(p+2)}\d^{p/(p+2)}=C'_{13} \varrho^{2/(p+2)}\d^{p/(p+2)},$$
where $C'_{13}=C_3^{p/(p+2)}(\xi+2\sqrt{2})^{p/(p+2)}.$

\noi
{\bf (ii).}  Let $q=0.$

Now repeating the arguments of (i) and using the fact $\Phi(\tilde{\b}_\d)=\xi\d^\nu$, we obtain
$$\|Tg-Tu^\d_{\b_\d}(0)\|\leq (\xi+2\sqrt{2})\d^\nu.$$ Also, it follows that $g-u_{\tilde{\b}_\d}(0)\in H_{p}$ and $\|g-u_{\tilde{\b}_\d}(0)\|_{H_{p}}\leq \varrho.$ Now, the proof follows by invoking the estimate in Theorem \ref{cond-stability}.
\epf
Now we are in a position to state the main result of this section, whose proof follows by combining the results in Theorems \ref{estimate-reg-noisy-aposteriori-quasi-bound} and \ref{estimate-reg-exact-aposteriori-quasi-bound}, respectively.
\bt\label{rate-quasi-bound-aposteriori}
Let $\d,\nu,u_{\b_\d}$ and $u_{\tilde{\b}_\d}$ be as in Lemma \ref{reg-parameter-bound-aposteriori-quasi-bound}. Let $g$ be the solution of \eqref{op-eq}. If $g\in S_{\varrho,p},$ for some $\varrho,p>0$ then the following holds:
\ben
\i[(i)] For $q\neq 0,\,0<p<q,$ we have $$\|g-u^\d_{\b_\d}(0)\|\leq C'_{13}\varrho^{2/(p+2)}\d^{p/(p+2)}+C_{14} \varrho^{2/(p+2)}\d^{p/(p+2)}.$$
\i[(ii)] For $q\neq 0,\,p\geq q,$ we have $$\|g-u^\d_{\b_\d}(0)\|\leq C'_{13}\varrho^{2/(p+2)}\d^{p/(p+2)}+C_{15}\varrho^{2/(q+2)}\d^{q/(q+2)}.$$
\i[(iii)] For $q=0,$ we have $$\|g-u^\d_{\tilde{\b}_\d}(0)\|\leq C'_{13} \varrho^{2/(p+2)}\d^{p\nu/(p+2)}+C_{16}\varrho\d^{1-\nu}.$$
\een
\et
\section{Fourier truncation method}\label{sec-FTM}
In this section we obtain regularized approximations by truncating the Fourier expansion of appropriate $L^2$ functions. The level at which the series expansion is truncated plays the role of regularization parameter. We shall provide both the apriori and aposteriori parameter choice strategies and derive the corresponding convergence rates of the approximation. We shall see that the rates obtained by this method are similar for both the apriori and aposteriori parameter choice. Moreover, we shall see that the rates obtained by this method does not possess the saturation effect which is in contrast to the rates obtained by quasi-boundary value method and its modifications, as obtained in the previous section. 

Throughout this section we assume the following: Let $\d>0,h,h^\d,f,f^\d$ be as in Theorem \ref{stability-initial-quasi-bound}, that is, $h,h^\d\in L^2(\O)$, $f,f^\d\in L^\infty(0,\t;L^2(\O))$ and 
\beq\label{noise-truncation}
\|h-h^\d\|^2+\theta\|f-f^\d\|^2_{L^\infty(0,\t;L^2(\O))}\leq \d^2,
\eeq where $\theta=\sum_{n=1}^\infty\frac{1}{\l_n^4}.$

 Let $\Upsilon_{h,f}^{\a,n}$, $\Upsilon_{h^\d,f^\d}^{\a,n}$ be as defined in \eqref{psifan} and $g\in L^2(\O)$ be the solution of \eqref{op-eq}.  Recall that solving \eqref{op-eq} is an ill-posed problem \w.r.t. the perturbations in the data $h$ and $f$. Thus, we aim to obtain stable approximations for the sought initial value $g$ by the Fourier truncation method. Towards this, for $N\in\N,$ we define
\beq\label{truncated-version}
\begin{cases}
g^N&:=\sum_{n=1}^N \frac{\Upsilon_{h,f}^{\a,n}}{\Ea(-\l_n^2\t^\a)}\f_n,\\
g^{\d,N}&:=\sum_{n=1}^N \frac{\Upsilon_{h^\d,f^\d}^{\a,n}}{\Ea(-\l_n^2\t^\a)}\f_n.
\end{cases}
\eeq
Here the natural number $N$ plays the role of regularization parameter. We shall choose $N$ in terms of $\d$ so that under certain source condition we can derive the error estimates in terms of the noise $\d.$ 
First, we obtain an estimate for $\|g-g^N\|.$ For $p,\varrho>0,$ recall the source set $S_{\varrho,p}$ as defined in \eqref{source-set}.
\bt\label{estimate-apriori-exact-truncation}
Let $g$ be the solution of \eqref{op-eq} and for $N\in\N,$ let $g^N$ be as defined in \eqref{truncated-version}. If $g\in S_{\varrho,p}$ for some $\varrho,p>0$ then we have$$\|g-g^N\|\leq \frac{\varrho}{\l_{N+1}^{p}}.$$
\et
\bpf
Since $g$ is a solution of \eqref{op-eq}, recall from \eqref{g-h-f-relation} that, we have $g=\sum_{n=1}^\infty g_n\f_n,$ where $g_n=\frac{\Upsilon_{h,f}^{\a,n}}{\Ea(-\l_n^2\t^\a)}.$ Thus, from the definition of $g^N,$ we have
\beqarray
\|g-g^N\|^2&=& \|\sum_{n=N+1}^\infty \frac{\Upsilon_{h,f}^{\a,n}}{\Ea(-\l_n^2\t^\a)}\f_n\|^2\\
&=& \sum_{n=N+1}^\infty \left(\frac{\Upsilon_{h,f}^{\a,n}}{\Ea(-\l_n^2\t^\a)}\right)^2\\
&=&\sum_{n=N+1}^\infty \frac{1}{\l_n^{2p}} \l_n^{2p}g_n^2\leq \frac{1}{\l_{N+1}^{2p}}\|g\|_{H_{p}}^2\\
&\leq &\frac{\varrho^2}{\l_{N+1}^{2p}},
\eeqarray
that is, $$\|g-g^N\|\leq \frac{\varrho}{\l_{N+1}^{p}}.$$
\epf
We now obtain bound for $\|g^N-g^{\d,N}\|.$
\bt\label{estimate-apriori-noisy-truncation}
For $\d>0$ and $N\in\N,$ let $g^N$ and $g^{\d,N}$ be as defined in \eqref{truncated-version} with $h,h^\d,f,f^\d$ satisfying the noise level as considered in \eqref{noise-truncation}. Then there exists a constant $C_{18}>0,$ depending on $\a,\t,\l_1,C_1,$ such that $$\|g^N-g^{\d,N}\|\leq C_{18}\d\l_N^2.$$
\et
\bpf
First we observe that from Theorem \ref{prop_Ea_upper_bound_lower_bound}, we have $$\frac{1}{\Ea(-\l_n^2\t^\a)}\leq \tilde{C}_{18}\l_n^2,$$ where $\tilde{C}_{18}=\frac{\Gamma(1-\a)(1+\t^\a\l_1^2)}{C_1\l_1^2}.$ Now from the definition of $g^N$ and $g^{\d,N},$ we have
\beqarray
\|g^N-g^{\d,N}\|^2&=&\|\sum_{n=1}^N \frac{\Upsilon_{h,f}^{\a,n}-\Upsilon_{h^\d,f^\d}^{\a,n}}{\Ea(-\l_n^2\t^\a)}\f_n\|^2=\sum_{n=1}^N \left(\frac{\Upsilon_{h,f}^{\a,n}-\Upsilon_{h^\d,f^\d}^{\a,n}}{\Ea(-\l_n^2\t^\a)}\right)^2\\
&\leq &  {\tilde{C}_{18}}^2\sum_{n=1}^N \l_n^4 (\Upsilon_{h,f}^{\a,n}-\Upsilon_{h^\d,f^\d}^{\a,n})^2\\
&\leq & {\tilde{C}_{18}}^2 \l_N^4\sum_{n=1}^\infty (\Upsilon_{h,f}^{\a,n}-\Upsilon_{h^\d,f^\d}^{\a,n})^2\\
&\leq & 2 \tilde{C_{18}}^2 \l_N^4 \d^2,
\eeqarray
where $C_{18}=\sqrt{2}\tilde{C}_{18}.$ Therefore,
$$\|g^N-g^{\d,N}\|\leq C_{18}\d \l_N^2.$$
\epf
\brem
Observe that Theorem \ref{estimate-apriori-noisy-truncation} shows that $g^N$'s  are stable \w.r.t. the perturbations in the data $h$ and $f.$
\erem
We are now in a position to state and prove one of the main result of this section, that is, the apriori parameter choice strategy and the corresponding error estimates. From now onwards we shall use the notation $[[r]]$ to denote the greatest integer not exceeding $r,$ where $r\in\R.$
\bt\label{apriori-rates-truncation}
Let $e_{1}$ and $e_2$ be as in \eqref{eig-value-bound}. For $\d>0$ and $N\in\N,$ let $h,h^\d,f,f^\d,g^{\d,N}$ be as in Theorem \ref{estimate-apriori-noisy-truncation}. For $\varrho,p>0,$ let $g$ be as in Theorem \ref{estimate-apriori-exact-truncation}. Let $C_{18}$ be as in Theorem \ref{estimate-apriori-noisy-truncation}. Then for the choice $$N_\d=\left[\left[\left(\frac{1}{C_{18}e_1^{p}e_2^2}\frac{\varrho}{\d}\right)^{d/(2p+4)}\right]\right],$$ there exists a constant $C_{19}>0,$ depending on $\a,\t,\l_1,C_1,e_1,e_2,$ such that $$\|g-g^{\d,N_\d}\|\leq C_{19}\varrho^{2/(p+2)}\d^{p/(p+2)}.$$
\et
\bpf
Combining the estimates obtained in Theorems \ref{estimate-apriori-exact-truncation} and \ref{estimate-apriori-noisy-truncation}, we have 
\beqarray
\|g-g^{\d,N}\|&\leq & \frac{\varrho}{\l_{N+1}^{p}}+C_{18}\d\l_N^2\leq \frac{\varrho}{\l_N^{p}}+C_{18}\d\l_N^2\\
&\underbrace{\leq}_{\eqref{eig-value-bound}} &\frac{\varrho}{e_1^{p}N^{2p/d}}+C_{18}\d e_2^2 N^{4/\d}.
\eeqarray
For $r>0,$ let $$S(r)=\frac{\varrho}{e_1^{p}r^{2p/d}}+C_{18}e_2^2\d r^{4/d}.$$ Then it is easy to verify that $S(r)$ attains its minimum at $r_\d$ satisfying $$\frac{\varrho}{e_1^{p}r_\d^{2p/d}}=C_{18}e_2^2\d r_\d^{4/d}.$$ That is, we have, $$r_\d=\left(\frac{1}{C_{18}e_1^{p}e_2^2}\frac{\varrho}{\d}\right)^{d/(2p+4)}.$$
We choose the regularization parameter $N_\d\in\N$ as $N_\d=\left[\left[r_\d\right]\right],$ that is,  $$N_\d=\left[\left[\left(\frac{1}{C_{18}e_1^{p}e_2^2}\frac{\varrho}{\d}\right)^{d/(2p+4)}\right]\right].$$
Then with this choice of $N_\d,$ we have,
$$\|g-g^{\d,N}\|\leq C_{18}e_2^2\,\d\left(\frac{1}{C_{18}e_1^{p}e_2^2}\,\frac{\varrho}{\d}\right)^{2/(p+2)}=C_{19}\varrho^{2/(p+2)}\d^{p/(p+2)},$$
where $C_{19}=\left(\frac{C_{18}e_2^2}{e_1^2}\right)^{p/(p+2)}.$ 

This completes the proof.
\epf
In Theorem \ref{apriori-rates-truncation} we have provided a parameter choice strategy and obtained the corresponding convergence rate. However, we observe that the choice of the regularization parameter $N_\d$ depends on $p,$ which is related to the smoothness of the unknown $g.$ Therefore, the parameter choice strategy will work only if we know apriori the smoothness of $g.$ But, in practice, it is almost always impossible to know the smoothness apriori. Thus, to overcome this situation, we provide the aposteriori parameter choice that only requires knowledge of the noise $\d$ and the observed data $h^\d,f^\d.$ 
\subsection{Error estimates with aposteriori parameter choice}
Recall that from Proposition \ref{combined-exact-data}, it follows that for $h\in L^2(\O)$ and $f\in L^\infty(0,\t;L^2(\O)),$ we have that $\Upsilon_{h,f}=\sum_{n=1}^\infty \Upsilon_{h,f}^{\a,n}\,\f_n\in L^2(\O),$ where $\Upsilon_{h,f}^{\a,n}$ is as defined in \eqref{psifan}. For $N\in \N,$ we let $$\Upsilon^N_{h,f}=\sum_{n=1}^N \Upsilon_{h,f}^{\a,n}\,\f_n,\q\text{and}\q\Upsilon^N_{h^\d,f^\d}=\sum_{n=1}^N \Upsilon_{h^\d,f^\d}^{\a,n}\,\f_n.$$ We define $$\zeta(N)=\|\Upsilon_{h^\d,f^\d}-\Upsilon^N_{h^\d,f^\d}\|,\q N\in\N.$$ Then it is easily seen that $\zeta(N)\to 0$ as $N\to \infty.$ 

Let $\mu>\sqrt{2}.$ Then, for $\d>0$ sufficiently small, there exists $N\in \N$ such that
\beq\label{discrepancy-principle-Fourier}
\|\Upsilon_{h^\d,f^\d}-\Upsilon^N_{h^\d,f^\d}\|\leq \mu \d<\|\Upsilon_{h^\d,f^\d}-\Upsilon^{N-1}_{h^\d,f^\d}\|.\eeq
Let $N_\d$ be the first natural number that satisfies \eqref{discrepancy-principle-Fourier}. This $N_\d$ will be the regularization parameter that we are looking for. Note that, clearly the choice of such  $N_\d$ depends only on the known quantities, $\d,h^\d,f^\d.$

We now obtain a bound for $N_\d$ in terms of $\d.$ 
\bl\label{estimate-reg-parameter-in-terms-of-delta-truncation}
Let $e_1$ be as in \eqref{eig-value-bound} and $h,f$ be the exact data. For $\d>0,$ let the noisy data $h^\d$ and $f^\d$ be such that it satisfies \eqref{noise-truncation}. Let $\mu>\sqrt{2}$ be fixed, $N_\d$ be the smallest natural number satisfying \eqref{discrepancy-principle-Fourier} and $C_2$ be as in Theorem \ref{prop_Ea_upper_bound_lower_bound}. Let $g$ be the solution of \eqref{op-eq}. If $g\in S_{\varrho,p}$ for some $p,\varrho>0,$ then there exists a constant $C_{20}>0,$ depending on $\a,\t,C_2,$ such that $$N_\d\leq \left(\frac{C_{20}}{(\mu-\sqrt{2})e_1^{(p+2)}}\,\frac{\varrho}{\d}\right)^{d/(2p+4)}.$$
\el
\bpf
First we note that $$\Upsilon_{h,f}-\Upsilon^{N_\d-1}_{h,f}=\sum_{n=N_\d}^\infty \Upsilon_{h,f}^{\a,n}\,\f_n=\sum_{n=N_\d}^\infty \Ea(-\l_n^2\t^\a)g_n\f_n.$$ Thus,
\beqarray
\|\Upsilon_{h,f}-\Upsilon^{N_\d-1}_{h,f}\|^2 &=& \sum_{n=N_\d}^\infty |\Ea(-\l_n^2\t^\a)|^2g_n^2\\
&\underbrace{\leq }_{Theorem \ref{prop_Ea_upper_bound_lower_bound}}& \sum_{n=N_\d}^\infty \frac{C_2^2}{(\Gamma(1-\a))^2\t^{2\a}\l_n^4}g_n^2=C_{20}^2\sum_{n=N_\d}^\infty \frac{g_n^2}{\l_n^4}\\
&=&C_{20}^2\sum_{n=N_\d}^\infty \frac{\l_n^{2p}g_n^2}{\l_n^{2p+4}}\leq C_{20}^2\frac{1}{\l_{N_\d}^{2p+4}}\sum_{n=1}^\infty \l_n^{2p}g_n^2\\
&\leq &C_{20}^2\frac{\varrho^2}{\l_{N_\d}^{2p+4}},
\eeqarray
where $C_{20}=\frac{C_2}{\Gamma(1-\a)\t^\a}.$ Thus, $$\|\Upsilon_{h,f}-\Upsilon^{N_\d-1}_{h,f}\|\leq C_{20}\frac{\varrho}{\l_{N_\d}^{p+2}}\leq \frac{C_{20}}{e_1^{p+2}}\frac{\varrho}{N_\d^{(2p+4)/d}}.$$
Now, observe that $$\|(\Upsilon_{h,f}-\Upsilon_{h^\d,f^\d})-(\Upsilon^{N_\d-1}_{h,f}-\Upsilon^{N_\d-1}_{h^\d,f^\d})\|^2=\sum_{n=N_\d}^\infty |\Upsilon_{h,f}^{\a,n}-\Upsilon_{h^\d,f^\d}^{\a,n}|^2\leq 2\d^2.$$
Therefore, we have
\beqarray
\|\Upsilon_{h,f}-\Upsilon^{N_\d-1}_{h,f}\|&=&\|(\Upsilon_{h,f}-\Upsilon_{h^\d,f^\d})-(\Upsilon^{N_\d-1}_{h,f}-\Upsilon^{N_\d-1}_{h^\d,f^\d})+(\Upsilon_{h^\d,f^\d}-\Upsilon^{N_\d-1}_{h^\d,f^\d})\|\\
&\geq & \|\Upsilon_{h^\d,f^\d}-\Upsilon^{N_\d-1}_{h^\d,f^\d}\|-\|(\Upsilon_{h,f}-\Upsilon_{h^\d,f^\d})-(\Upsilon^{N_\d-1}_{h,f}-\Upsilon^{N_\d-1}_{h^\d,f^\d})\|\\
&\geq & (\mu-\sqrt{2})\d.
\eeqarray
Thus,
$$(\mu-\sqrt{2})\d\leq \frac{C_{20}}{e_1^{p+2}}\frac{\varrho}{N_\d^{(2p+4)/d}},$$
and hence $$N_\d\leq \left(\frac{C_{20}}{(\mu-\sqrt{2})e_1^{(p+2)}}\,\frac{\varrho}{\d}\right)^{d/(2p+4)}.$$
\epf

We now obtain estimate for $\|g-g^{N_\d}\|.$
\bt\label{estimate-exact-reg-aposteriori-truncation}
Let $g$ be the solution of \eqref{op-eq} and $g\in S_{\varrho,p}$ for some $\varrho,p>0.$ Let $\mu,\d, N_\d$ be as in Lemma \ref{estimate-reg-parameter-in-terms-of-delta-truncation} and $g^{N_\d}$ be as defined in \eqref{truncated-version}. Let $C_1$ be as in Theorem \ref{prop_Ea_upper_bound_lower_bound}. Then there exists a constant $C_{21}>0,$ depending on $\a,\t,\l_1,C_1,p,\mu,$ such that $$\|g-g^{N_\d}\|\leq C_{21}\varrho^{2/(p+2)}\d^{p/(p+2)}.$$
\et
\bpf
Since $g$ is a solution of \eqref{op-eq}, recall that $g=\sum_{n=1}^\infty g_n\f_n,$ where $g_n$ is as given in \eqref{g-h-f-relation}, that is, $g_n=\frac{\Upsilon_{h,f}^{\a,n}}{\Ea(-\l_n^2\t^\a)}.$
Thus, from the definition of $g^{N_\d},$ we have
\beqarray
\|g-g^{N_\d}\|^2 &=& \sum_{n=N_\d+1}^\infty \left(\frac{\Upsilon_{h,f}^{\a,n}}{\Ea(-\l_n^2\t^\a)}\right)^2\\
&=& \sum_{n=N_\d+1}^\infty \frac{(\Upsilon_{h,f}^{\a,n})^{4/(p+2)}}{(\Ea(-\l_n^2\t^\a))^2} (\Upsilon_{h,f}^{\a,n})^{2p/(p+2)}\\
&\leq & \left(\sum_{n=N_\d+1}^\infty\frac{(\Upsilon_{h,f}^{\a,n})^2}{(\Ea(-\l_n^2\t^\a))^{p+2}}\right)^{2/(p+2)}\left(\sum_{n=N_\d+1}^\infty(\Upsilon_{h,f}^{\a,n})^2\right)^{p/(p+2)}\\
&\leq & \left(\sum_{n=1}^\infty \frac{g_n^2}{(\Ea(-\l_n^2\t^\a))^{p}}\right)^{2/(p+2)}\left(\sum_{n=N_\d+1}^\infty (\Upsilon_{h,f}^{\a,n})^2\right)^{p/(p+2)}\\
&\leq &\left({\tilde{C}_{18}}^{p}\sum_{n=1}^\infty \l_n^{2p} g_n^2\right)^{2/(p+2)}\left(\sum_{n=N_\d+1}^\infty (\Upsilon_{h,f}^{\a,n})^2\right)^{p/(p+2)}\\
&\leq & {\tilde{C}_{18}}^{2p/(p+2)} \varrho^{4/(p+2)}\left(\sum_{n=N_\d+1}^\infty (\Upsilon_{h,f}^{\a,n})^2\right)^{p/(p+2)},
\eeqarray
where $\tilde{C}_{18}=\frac{\Gamma(1-\a)(1+\t^\a\l_1^2)}{C_1\l_1^2}.$

Now, recall that $$\|\Upsilon_{h^\d,f^\d}-\Upsilon^{N_\d}_{h^\d,f^\d}\|\leq \mu\d$$ and $$\|(\Upsilon_{h,f}-\Upsilon_{h^\d,f^\d})-(\Upsilon^{N_\d}_{h,f}-\Upsilon^{N_\d}_{h^\d,f^\d})\|^2=\sum_{n=N_\d+1}^\infty |\Upsilon_{h,f}^{\a,n}-\Upsilon_{h^\d,f^\d}^{\a,n}|^2\leq 2\d^2.$$
Thus,
\beqarray
\sum_{n=N_\d+1}^\infty (\Upsilon_{h,f}^{\a,n})^2 &=& \|\Upsilon_{h,f}-\Upsilon^{N_\d}_{h,f}\|^2\\
&=&\|(\Upsilon_{h^\d,f^\d}-\Upsilon^{N_\d}_{h^\d,f^\d})+(\Upsilon_{h,f}-\Upsilon_{h^\d,f^\d})-(\Upsilon^{N_\d}_{h,f}-\Upsilon^{N_\d}_{h^\d,f^\d})\|^2\\
&\leq & (\mu+\sqrt{2})^2\d^2.
\eeqarray
Thus, we obtain
$$\|g-g^{N_\d}\|^2\leq {\tilde{C}_{18}}^{2p/(p+2)}\varrho^{4/(p+2)}(\mu+\sqrt{2})^{2p/(p+2)}\d^{2p/(p+2)},$$
that is,
$$\|g-g^{N_\d}\|\leq C_{21}\varrho^{2/(p+2)}\d^{p/(p+2)},$$
where $C_{21}={\tilde{C}_{18}}^{p/(p+2)}(\mu+\sqrt{2})^{p/(p+2)}.$
\epf
We are now in a position to state and prove the main result of this section, in which we obtain the rates of convergence corresponding to the aposteriori choice of the parameter $N_\d$.
\bt\label{aposteriori-rates-truncation}
Let $\mu>\sqrt{2}, \d>0,$ and $N_\d$ be as in Lemma \ref{estimate-reg-parameter-in-terms-of-delta-truncation}. Let $g^{\d,N_\d}$ be as in \eqref{truncated-version} and $e_1,e_2$ be as in \eqref{eig-value-bound}. Let $g$ be the solution of \eqref{op-eq}. If $g\in S_{\varrho,p},$ for some $\varrho,p>0,$ then there exists a constant $C_{23}>0,$ depending on $\mu,\a,\t,\l_1,e_1,e_2,p,$ such that $$\|g-g^{\d,N_\d}\|\leq C_{23}\varrho^{2/(p+2)}\d^{p/(p+2)}.$$
\et
\bpf
From the estimate obtained in Theorem \ref{estimate-apriori-noisy-truncation}, we have $$\|g^{N_\d}-g^{\d,N_\d}\|\leq C_{18}\d\l_{N_\d}^2\leq C_{22}\d N_\d^{4/d},$$ where $C_{22}=e_2^2 C_{18}$ and $e_2$ is as in \eqref{eig-value-bound}.

Thus, from Theorem \ref{estimate-exact-reg-aposteriori-truncation}, we have
\beqarray
\|g-g^{\d,N_\d}\|&\leq &\|g-g^{N_\d}\|+\|g^{N_\d}-g^{\d,N_\d}\|\\
&\leq & C_{21}\varrho^{2/(p+2)}\d^{p/(p+2)}+C_{22}\,\d N_\d^{4/d}\\
&\leq & C_{21}\varrho^{2/(p+2)}\d^{p/(p+2)}+ C_{22}\,\d\left(\frac{C_{20}}{(\mu-\sqrt{2})e_1^{p+2}}\right)^{2/(p+2)}\left(\frac{\varrho}{\d}\right)^{2/(p+2)}\\
&=& C_{21}\varrho^{2/(p+2)}\d^{p/(p+2)}+C_{22}\left(\frac{C_{20}}{(\mu-\sqrt{2})e_1^{p+2}}\right)^{2/(p+2)}\varrho^{2/(p+2)}\d^{p/(p+2)}\\
&=&C_{23} \varrho^{2/(p+2)}\d^{p/(p+2)},
\eeqarray
where $C_{23}=\max\big\{C_{21},C_{22}\big(\frac{C_{20}}{(\mu-\sqrt{2})e_1^{p+2}}\big)^{2/(p+2)}\big\}.$
\epf
\brem
Clearly the rates obtained in Theorems \ref{apriori-rates-truncation} and \ref{aposteriori-rates-truncation}, for the apriori and aposteriori case, respectively, are of the same order, that is, $O(\d^{p/(p+2)})$ for all $p>0.$ Since $\d^{p/(p+2)}\to \d$ as $p\to \infty,$ it follows that the rates for both apriori and aposteriori, obtained by the Fourier truncation method does not possesses the saturation effect.
\erem

\brem{\bf (Optimality)}\label{optimality}
In Section \ref{sec-prelim}, we have already observed that the inverse problem of retrieving the initial value $g=u(\cdot,0)$ from the knowledge of the final value $h=u(\cdot,\t)$ and the source function $f$ is same as solving the linear operator equation $$Tg=\Upsilon_{h,f},$$
where $\Upsilon_{h,f}=\sum_{n=1}^\infty \Upsilon^{\a,n}_{h,f}\f_n,$ $\Upsilon_{h,f}^{\a,n}$ is as defined in \eqref{psifan}, and the compact linear self-adjoint operator $T:L^2(\O)\to L^2(\O)$ is defined by $$T\phi=\sum_{n=1}^\infty \k_n\<\phi,\f_n\>\f_n,\q\phi\in L^2(\O),$$ where $\k_n=\Ea(-\l_n^2\t^\a).$ 

From Theorem \ref{prop_Ea_upper_bound_lower_bound}, we have $\k_n>0$ for all $n\in\N,$ and $$\frac{\tilde{C}_1}{\l_n^2}\leq \k_n\leq \frac{\tilde{C}_2}{\l_n^2},\q \forall n\in\N,$$
where $\tilde{C}_1=\frac{C_1\l_1^2}{(1+\l_1^2\t^\a)\Gamma(1-\a)}$ and $\tilde{C}_2=\frac{C_2}{\t^\a\Gamma(1-\a)}.$ This shows that, $T$ is also a self-adjoint positive operator. Now following the arguments in \cite[Section 5]{mondal_nair_2024}, we can arrive at the conclusion that for the source set $S_{\varrho,p}$ and the operator $T$ as considered above, the convergence rate $O(\d^{p/(p+2)})$ is order optimal. 

Thus, from Theorems \ref{apriori-rates-truncation} and \ref{aposteriori-rates-truncation}, it follows that for the source set $S_{\varrho,p},$ the convergence rates obtained by FTM are order optimal for both the apriori and aposteriori parameter choice strategies for all $p>0$. Also, from Theorems \ref{quasi-bound-apriori-rate} and \ref{rate-quasi-bound-aposteriori}, we observe that for MQBVM: if $0<p\leq q+2$ then the rate obtained for the apriori case  is order optimal, whereas if $0<p\leq q$ then the rate obtained for the apsoteriori case is order optimal. 
\erem
\section{Conclusion}
We have initiated a study on the backward problem of retrieving initial value from final value for a time-fractional fourth order parabolic PDE. Since the considered problem is ill-posed, we have obtained stable approximations by employing quasi-boundary value method, its modifications and the Fourier truncation method. Under some Sobolev smoothness assumption on the sought initial value, we have obtained convergence rates for all these methods under both apriori and aposteriori parameter choice strategies. We observed that for certain cases all the considered methods produces same convergence rates in both the apriori and aposteriori cases. However, we observed that the rates obtained by QBVM and MQBVM possesses the saturation effect but the rates obtained by FTM is free of such effect. Moreover, we observed that the rates obtained by FTM is order optimal for all values of the smoothness index $p$.

\vspace{1cm}
\noi
{\bf Data availability} Data sharing is not applicable to this article.

\vspace{.5cm}
\noi
{\bf Declarations} 

\vspace{.3cm}
\noi
{\bf Conflict of interest} The author declares that there is no conflict of interest.

\vspace{.5cm}
\noi
{\bf Acknowledgements.} The author is supported by the postdoctoral fellowship of the TIFR Centre for Applicable Mathematics, Bangalore.

\end{document}